\documentclass[leqno,english]{smfart}

\usepackage[leqno]{amsmath}
\usepackage{amssymb, amsthm, enumerate, latexsym, mathrsfs, mdwlist, stmaryrd}
\usepackage{ae,aecompl}
\usepackage[all]{xy}

\theoremstyle{remark}
\newtheorem{nul}{}[section]

\newtheorem*{rem*}{Remark}

\theoremstyle{definition}
\newtheorem{dfn}[nul]{Definition}

\newtheorem{ntn}[nul]{Notation}
\newtheorem{exm}[nul]{Example}

\newtheorem*{dfn*}{Definition}
\newtheorem*{axm*}{Axiom}
\newtheorem*{ntn*}{Notation}
\newtheorem*{exm*}{Example}
\newtheorem*{exr*}{Exercise}
\newtheorem*{int*}{Intuition}
\newtheorem*{qst*}{Question}

\theoremstyle{plain}
\newtheorem{sch}[nul]{Scholium}
\newtheorem{thm}[nul]{Theorem}
\newtheorem{prp}[nul]{Proposition}
\newtheorem{cor}[nul]{Corollary}
\newtheorem{lem}[nul]{Lemma}

\newtheorem*{thm*}{Theorem}
\newtheorem*{prp*}{Proposition}
\newtheorem*{cor*}{Corollary}
\newtheorem*{lem*}{Lemma}
\newtheorem*{cnj*}{Conjecture}

\numberwithin{equation}{nul}

\swapnumbers

\DeclareMathOperator{\cell}{cell}
\DeclareMathOperator{\cof}{cof}

\DeclareMathOperator{\const}{const}
\DeclareMathOperator{\colim}{colim}

\DeclareMathOperator{\diag}{diag}

\DeclareMathOperator{\fib}{fib}

\DeclareMathOperator{\Ho}{Ho}
\DeclareMathOperator{\hocolim}{hocolim}
\DeclareMathOperator{\holim}{holim}
\DeclareMathOperator{\id}{id}

\DeclareMathOperator{\inj}{inj}

\DeclareMathOperator{\Map}{Map}
\DeclareMathOperator{\mor}{mor}
\DeclareMathOperator{\Mor}{Mor}
\DeclareMathOperator{\MOR}{\mathbf{Mor}}

\DeclareMathOperator{\Obj}{Obj}
\DeclareMathOperator{\Post}{Post}

\DeclareMathOperator{\proj}{proj}

\DeclareMathOperator{\RMor}{\mathbf{R}Mor}
\DeclareMathOperator{\RMOR}{\mathbf{RMor}}

\DeclareMathOperator{\Sect}{Sect}

\newcommand{\Prod}{\prod}
\newcommand{\Coprod}{\coprod}

\newcommand{\coprd}{\amalg}

\newcommand{\CC}{\mathbf{C}}
\newcommand{\DD}{\mathbf{D}}
\newcommand{\EE}{\mathbf{E}}
\newcommand{\FF}{\mathbf{F}}
\newcommand{\GG}{\mathbf{G}}
\newcommand{\HH}{\mathbf{H}}

\newcommand{\LL}{\mathbf{L}}
\newcommand{\MM}{\mathbf{M}}
\newcommand{\NN}{\mathbf{N}}

\newcommand{\PP}{\mathbf{P}}

\newcommand{\RR}{\mathbf{R}}

\newcommand{\VV}{\mathbf{V}}

\newcommand{\XX}{\mathbf{X}}
\newcommand{\YY}{\mathbf{Y}}
\newcommand{\ZZ}{\mathbf{Z}}

\newcommand{\op}{\mathrm{op}}

\newcommand{\fromto}[2]{\xymatrix@1@C=18pt{{#1}\ar[r]&{#2}}}
\newcommand{\goesto}[2]{\xymatrix@1@C=18pt{{#1}\,\ar@{|->}[r]&{#2}}}

\begin{document}

\title[On (Enriched) Left Bousfield Localizations]{On (Enriched) Left Bousfield Localizations of Model Categories}
\author{Clark Barwick}
\address{Matematisk Institutt\\
Universitetet i Oslo\\
Boks 1053 Blindern\\
0316 Oslo\\
Norge}
\curraddr{School of Mathematics\\
Institute for Advanced Study\\
Einstein Drive\\
Princeton, NJ 08540-0631\\
USA}
\email{clarkbar@gmail.com}
\thanks{This work was supported by a research grant from the Yngre fremragende forskere, administered by J. Rognes at the Matematisk Institutt, Universitetet i Oslo.}

\begin{abstract} I verify the existence of left Bousfield localizations and of enriched left Bousfield localizations, and I prove a collection of useful technical results characterizing certain fibrations of (enriched) left Bousfield localizations. I also use such Bousfield localizations to construct a number of new model categories, including models for the homotopy limit of right Quillen presheaves, for Postnikov towers in model categories, and for presheaves valued in a symmetric monoidal model category satisfying a homotopy-coherent descent condition.
\end{abstract}

\maketitle
\thispagestyle{empty}

A class of maps $H$ in a model category $\MM$ specifies a class of $H$-local objects, which are those objects $X$ with the property that the morphism $\RMor_{\MM}(f,X)$ is a weak equivalence of simplicial sets for any $f\in H$. The left Bousfield localization of $\MM$ with respect to $H$ is a model for the homotopy theory of $H$-local objects. Similarly, if $\MM$ is enriched over a symmetric monoidal model category $\VV$, the class $H$ specifies a class of $(H/\VV)$-local objects, which are those objects $X$ with the property that the morphism $\RMOR_{\MM}^{\VV}(f,X)$ of derived mapping objects is a weak equivalence of $\VV$ for any $f\in H$. The $\VV$-enriched left Bousfield localization of $\MM$ is a model for the homotopy theory of $(H/\VV)$-local objects. 

The (enriched) left Bousfield localization is described as a new ``$H$-local'' model category structure on the underlying category of $\MM$. The $H$-local cofibrations of the (enriched) left Bousfield localization are precisely those of $\MM$, the $H$-local fibrant objects are the (enriched) $H$-local objects that are fibrant in $\MM$, and the $H$-local weak equivalences are those morphisms $f$ of $\MM$ such that $\RMor_{\MM}(f,X)$ (resp., $\RMOR_{\MM}^{\VV}(f,X)$) is a weak equivalence. This is enough to specify $H$-local fibrations, but it can be difficult to get explicit control over them. Luckily, is frequently possible to characterize some of the $H$-local fibrations as fibrations that are in addition homotopy pullbacks of fibrations between $H$-local fibrant objects.

The (enriched) Bousfield localization gives an effective way of constructing new model categories from old. In particular, one can use this to construct models for the homotopy limit of a right Quillen presheaf, for Postnikov towers in model categories, and for presheaves valued in a symmetric monoidal model category satisfying a homotopy-coherent descent condition.

The aim of this short note is to demonstrate the existence of (enriched) left Bousfield localizations, to give some useful technical results concerning the fibrations thereof, and to use this machinery to provide a number of interesting examples of model categories, including those described in the previous paragraph.

\subsection*{Plan} In the first section, I give a brief review of the general theory of combinatorial and tractable model categories. This material is all well-known. There one may find two familiar but important examples: model structures on diagram categories and model structures on section categories.

In the second section, I define the left Bousfield localization and give the well-known existence theorem due to Smith. Following this, I continue with a small collection of results that permit one to cope with the fact that left Bousfield localization ruins right properness, as well as a characterization of a certain class of $H$-local fibrations. I finish the section with five simple applications of the technique of left Bousfield localization: Dugger's presentation theorem, the existence of homotopy images, the existence of resolutions of model categories, the construction of homotopy limits of diagrams of model categories, and the existence of Postnikov towers for simplicial model categories.

In the final section, I review the notions of symmetric monoidal and enriched model categories. I describe the enriched left Bousfield localization and prove an existence theorem. I then give two applications: first, the existence of enriched Postnikov towers, and second, the existence of local model structures on presheaves valued in symmetric monoidal model categories.

Thanks to J. Bergner, P. A. {\O}stv{\ae}r, and B. Toën for persistent encouragement and hours of interesting discussion. Thanks also to J. Rosick{\'y} for pointing out a careless omission. Thanks especially to M. Spitzweck for a profound and lasting impact on my work; were it not for his insights and questions, there would be nothing for me to report here or anywhere else.

\setcounter{tocdepth}{2}
\tableofcontents

\section{Tractable model categories}

\subsection*{Combinatorial and tractable model categories} Combinatorial model categories are those whose homotopy theory is controlled by the homotopy theory of a small subcategory of presentable objects. A variety of algebraic applications require that the sets of (trivial) cofibrations can be generated (as a saturated set) by a given small set of (trivial) cofibrations with cofibrant domain. This leads to the notion of tractable model categories. Many of the results have satisfactory proofs in print; the first section of \cite{MR1780498} in particular is a very nice reference.\footnote{I would like to thank M. Spitzweck for suggesting this paper; this exposition has benefitted greatly from his recommendation.}

\begin{nul} Suppose here $\XX$ a universe.
\end{nul}

\begin{ntn}\label{ntn:lambdaprescoffib} Suppose $\CC$ a model $\XX$-category.
\begin{enumerate}[(\ref{ntn:lambdaprescoffib}.1)]
\item Write $w\CC$ (respectively, $\cof\CC$, $\fib\CC$) for the lluf subcategories comprised of weak equivalences (resp., cofibrations, fibrations).
\item Write $\CC_c$ (respectively, $\CC_f$) for the full subcategories of cofibrant (resp., fibrant) objects.
\item For any $\XX$-small regular cardinal $\lambda$, write $C_{\lambda}$ for the full subcategory of $C$ comprised of $\lambda$-presentable objects.
\end{enumerate}
\end{ntn}

\begin{dfn}\label{dfn:cofgensemi} Suppose $\CC$ a model $\XX$-category. Suppose, in addition, that $\lambda$ is a regular $\XX$-small cardinal.
\begin{enumerate}[(\ref{dfn:cofgensemi}.1)]
\item\label{item:cofgensemi}One says that $\CC$ is \emph{$\lambda$-combinatorial} (respectively, \emph{$\lambda$-tractable}) if its underlying $\XX$-category is locally $\lambda$-presentable, and if there exist $\XX$-small sets $I$ and $J$ of morphisms of $\CC_{\lambda}$ (resp., of $\CC_{\lambda}\cap\CC_c$) such that the following hold.
\begin{enumerate}[(\ref{dfn:cofgensemi}.\ref{item:cofgensemi}.1)]
\item A morphism satisfies the right lifting property with respect to $I$ if and only if it is a trivial fibration.
\item A morphism satisfies the right lifting property with respect to $J$ if and only if it is a fibration.
\end{enumerate}
\item An $\XX$-small full subcategory $\CC_0$ of $\CC$ is \emph{homotopy $\lambda$-generating} if every object of $\CC$ is weakly equivalent to a $\lambda$-filtered homotopy colimit of objects of $\CC_0$.
\item One says that $\CC$ is $\XX$-combinatorial (respectively, \emph{$\XX$-tractable}) just in case there exists a regular $\XX$-small cardinal $\kappa$ for which it is $\kappa$-combinatorial (resp., $\kappa$-tractable).
\item Likewise, an $\XX$-small full subcategory $\CC_0$ of $\CC$ is \emph{homotopy $\XX$-generating} if and only if for some regular $\XX$-small cardinal $\kappa$, $\CC_0$ is homotopy $\kappa$-generating.
\end{enumerate}
\end{dfn}

\begin{ntn} Suppose $C$ an $\XX$-category, and suppose $I$ an $\XX$-small set of morphisms of $C$. Denote by $\inj I$ the set of all morphisms with the right lifting property with respect to $I$, denote by $\cof I$ the set of all morphisms with the left lifting property with respect to $\inj I$, and denote by $\cell I$ the set of all transfinite compositions of pushouts of morphisms of $I$.
\end{ntn}

\begin{lem}[Transfinite small object argument, \protect{\cite[Proposition 1.3]{MR1780498}}]\label{lem_transfsoa} Suppose $\lambda$ a regular $\XX$-small cardinal, $C$ a locally $\lambda$-presentable $\XX$-category, and $I$ an $\XX$-small set of morphisms of $C_{\lambda}$.
\begin{enumerate}[(\ref{lem_transfsoa}.1)]
\item There is an accessible functorial factorization of every morphism $f$ as $p\circ i$, wherein $p\in\inj I$, and $i\in\cell I$.
\item A morphism $q\in\inj I$ if and only if it has the right lifting property with respect to all retracts of morphisms of $\cell I$.
\item A morphism $j$ is a retract of morphisms of $\cell I$ if and only $j\in\cof I$.
\end{enumerate}
\begin{proof} Suppose $\kappa$ a regular cardinal strictly greater than $\lambda$. For any morphism $f:\fromto{X}{Y}$, consider the set $(I/f)$ of squares
\begin{equation*}
\xymatrix@C=18pt@R=18pt{K\ar[d]_{i}\ar[r]&X\ar[d]^f\\
L\ar[r]&Y,}
\end{equation*}
where $i\in I$, and let $\fromto{K_{(I/f)}}{L_{(I/f)}}$ be the coproduct $\coprd_{i\in(I/f)}i$. Define a section $P$ of $d_1:\fromto{C^2}{C^1}$ by
\begin{equation*}
Pf:=[\xymatrix@1@C=18pt{X\ar[r]&X\sqcup^{L_{(I/f)}}K_{(I/f)}\ar[r]&Y}]
\end{equation*}
for any morphism $f:\fromto{X}{Y}$. For any regular cardinal $\alpha$, set $P^{\alpha}:=\colim_{\beta<\alpha}P^{\beta}$. This provides a functorial factorization $P^{\kappa}$ with the required properties.

The remaining parts follow from the existence of this factorization and the retract argument.
\end{proof}
\end{lem}

\begin{nul} J. Smith's insight is that the transfinite small object argument and the solution set condition on weak equivalences together provide a good recognition principle for combinatorial model categories. In effect, one requires only two-thirds of the data normally required to produce cofibrantly generated model structures.
\end{nul}

\begin{prp}[Smith, \protect{\cite[Theorem 1.7 and Propositions 1.15 and 1.19]{MR1780498}}]\label{lem_smith-beke} Suppose $C$ a locally $\XX$-presentable $\XX$-category, $W$ an accessibly embedded, accessible subcategory of $C^1$, and $I$ an $\XX$-small set of morphisms of $C$. Suppose in addition that the following conditions are satisfied.
\begin{enumerate}[(\ref{lem_smith-beke}.1)]
\item $W$ satisfies the two-out-of-three axiom.
\item The set $\inj I$ is contained in $W$.
\item The intersection $W\cap\cof I$ is closed under pushouts and transfinite composition.
\end{enumerate}
Then $C$ is a combinatorial model category with weak equivalences $W$, cofibrations $\cof I$, and fibrations $\inj(W\cap\cof I)$.
\begin{proof} The usual small object and retract arguments apply once one constructs an $\XX$-small set $J$ such that $\cof J=W\cap\cof I$. The following pair of lemmas complete the proof.
\end{proof}
\end{prp}

\begin{lem}[Smith, \protect{\cite[Lemma 1.8]{MR1780498}}]\label{lem:denseJ} Under the hypotheses of proposition \ref{lem_smith-beke}, suppose $J\subset W\cap\cof I$ a set such that any commutative square
\begin{equation*}
\xymatrix@C=18pt@R=18pt{
K\ar[d]\ar[r]&M\ar[d]\\
L\ar[r]&N
}
\end{equation*}
in which $[\fromto{K}{L}]\in I$ and $[\fromto{M}{N}]\in W$ can be factored as a commutative diagram
\begin{equation*}
\xymatrix@C=18pt@R=18pt{
K\ar[d]\ar[r]&M'\ar[d]\ar[r]&M\ar[d]\\
L\ar[r]&N'\ar[r]&N,
}
\end{equation*}
in which $[\fromto{M'}{N'}]\in J$. Then $\cof J=W\cap\cof I$.
\begin{proof} To prove this, one need only factor any element of $W$ as an element of $\cell J$ followed by an element of $\inj I$. The result then follows from the retract argument.

Suppose $\kappa$ an $\XX$-small regular cardinal such that every codomain of $I$ is $\kappa$-presentable. For any morphism $[f:\fromto{X}{Y}]\in W$, consider the set $(I/f)$ of squares
\begin{equation*}
\xymatrix@C=18pt@R=18pt{K\ar[d]_{i}\ar[r]&X\ar[d]^f\\
L\ar[r]&Y,}
\end{equation*}
where $i\in I$; for each such square choose an element $j_{(i,f)}\in J$ and a factorization
\begin{equation*}
\xymatrix@C=18pt@R=18pt{
K\ar[d]_{i}\ar[r]&M(i)\ar[d]^{j_{(i,f)}}\ar[r]&X\ar[d]^f\\
L\ar[r]&N(i)\ar[r]&Y,
}
\end{equation*}
and let $\fromto{M_{(I/f)}}{N_{(I/f)}}$ be the coproduct $\coprd_{i\in(I/f)}j_{(i,f)}$. Define an endofunctor $Q$ of $(W/Y)$ by
\begin{equation*}
Qf:=[\fromto{X\sqcup^{N_{(I/f)}}M_{(I/f)}}{Y}]
\end{equation*}
for any morphism $[f:\fromto{X}{Y}]\in W$. For any regular cardinal $\alpha$, set $Q^{\alpha}:=\colim_{\beta<\alpha}Q^{\beta}$. This provides, for any morphism $[f:\fromto{X}{Y}]\in W$, a functorial factorization
\begin{equation*}
\xymatrix@1@C=18pt{X\ar[r]&Q^{\kappa}f\ar[r]&Y}
\end{equation*}
with the desired properties.
\end{proof}
\end{lem}

\begin{lem}[Smith, \protect{\cite[Lemma 1.9]{MR1780498}}] Under the hypotheses of proposition \ref{lem_smith-beke}, an $\XX$-small set $J$ satisfying the conditions of lemma \ref{lem:denseJ} can be found.
\begin{proof} Since $W$ is an accessibly embedded accessible subcategory of $C^1$, it follows that for any morphism $[i:\fromto{K}{L}]\in I$, there exists an $\XX$-small subset $W(i)\subset W$ such that for any commutative square
\begin{equation*}
\xymatrix@C=18pt@R=18pt{
K\ar[d]_i\ar[r]&M\ar[d]\\
L\ar[r]&N
}
\end{equation*}
in which $[\fromto{M}{N}]\in W$, there exist a morphism $[\fromto{P}{Q}]\in W(i)$ and a commutative diagram
\begin{equation*}
\xymatrix@C=18pt@R=18pt{
K\ar[d]\ar[r]&P\ar[d]\ar[r]&M\ar[d]\\
L\ar[r]&Q\ar[r]&N.
}
\end{equation*}
It thus suffices to find, for every square of the type on the left, an element of $W\cap\cof I$ factoring it.

For every $[i:\fromto{K}{L}]\in I$, every $[w:\fromto{P}{Q}]\in W(i)$, and every commutative square
\begin{equation*}
\xymatrix@C=18pt@R=18pt{
K\ar[d]_i\ar[r]&P\ar[d]\\
L\ar[r]&Q,
}
\end{equation*}
factor the morphism $\fromto{L\sqcup^KP}{Q}$ as an element of $[\fromto{L\sqcup^KP}{R}]\in\cell I$ followed by an element of $[\fromto{R}{Q}]\in\inj I$; this yields a commutative diagram
\begin{equation*}
\xymatrix@C=18pt@R=18pt{
K\ar[d]\ar[r]&P\ar[d]\ar@{=}[r]&P\ar[d]\\
L\ar[r]&R\ar[r]&Q
}
\end{equation*}
factoring the original square, in which $[\fromto{P}{R}]\in W\cap\cof I$.
\end{proof}
\end{lem}

\begin{prp}[Smith, \protect{\cite[Propositions 7.1--3]{MR1870516}}]\label{lem_goodpropsofcombin} Suppose $\CC$ an $\XX$-combinatorial left model $\XX$-category. For any sufficient large $\XX$-small regular cardinal $\kappa$, the following hold.
\begin{enumerate}[(\ref{lem_goodpropsofcombin}.1)]
\item\label{item_kappacoftrivfib} There exists a $\kappa$-accessible functorial factorization $\fromto{\CC^1}{\CC^2}$ of each morphism into a cofibration followed by a trivial fibration.
\item\label{item_kappatrivcoffib} There exists a $\kappa$-accessible functorial factorization $\fromto{\CC^1}{\CC^2}$ of each morphism into a trivial cofibration followed by a fibration.
\item\label{item_kappacof} There exists a $\kappa$-accessible cofibrant replacement functor.
\item\label{item_kappafib} There exists a $\kappa$-accessible fibrant replacement functor.
\item\label{item_kappacolimwes} Arbitrary $\kappa$-filtered colimits preserve weak equivalences.
\item\label{item_kappacolimhocolim} Arbitrary $\kappa$-filtered colimits in $\CC$ are homotopy colimits.
\item\label{item_weaccessible} The set of weak equivalences $w\CC$ form a $\kappa$-accessibly embedded, $\kappa$-accessible subcategory of $\CC^1$.
\end{enumerate}
\begin{proof} Observe that (\ref{lem_goodpropsofcombin}.\ref{item_kappacoftrivfib}) and (\ref{lem_goodpropsofcombin}.\ref{item_kappatrivcoffib}) (and therefore (\ref{lem_goodpropsofcombin}.\ref{item_kappacof}) and (\ref{lem_goodpropsofcombin}.\ref{item_kappafib}) as well) follow directly from the transfinite small object argument \ref{lem_transfsoa}.

To verify (\ref{lem_goodpropsofcombin}.\ref{item_kappacolimwes}) --- and therefore (\ref{lem_goodpropsofcombin}.\ref{item_kappacolimhocolim}) ---, fix $\kappa$, an $\XX$-small regular cardinal for which: (a) there are $\kappa$-accessible functorial factorizations of each kind, (b) there is an $\XX$-small set $I$ of generating cofibrations with $\kappa$-presentable domains and codomains, and (c) the full subcategory of $I$-tuples of surjective morphisms is a $\kappa$-access9ibly embedded, $\kappa$-accessible subcategory of $\mathbf{Set}_{\XX}^{(I\cdot 1}$. Suppose $A$ an $\XX$-small $\kappa$-filtered category, and suppose $\fromto{F}{G}$ an objectwise weak equivalence in $\CC^A$. The $\kappa$-accessible functorial factorizations in $\CC$ permit one to give a $\kappa$-accessible factorization of $\fromto{F}{G}$ into an objectwise trivial cofibration $\fromto{F}{H}$ followed by an objectwise fibration (which is therefore an objectwise trivial fibration) $\fromto{H}{G}$. Hence the morphism $\fromto{\colim H}{\colim G}$ is a fibration, and it remains only to show that it is also a trivial fibration; for this one need only show that for any morphism $f:\fromto{K}{L}$ and element of $I$ and any diagram
\begin{equation*}
\xymatrix@C=18pt@R=18pt{K\ar[d]\ar[r]&\colim F\ar[d]\\
L\ar[r]&\colim G,}
\end{equation*}
a lift $\fromto{L}{\colim F}$ exists. This follows from the $\kappa$-presentability of $K$ and $L$.

To verify (\ref{lem_goodpropsofcombin}.\ref{item_weaccessible}), let us note that it follows from the existence of a $\kappa$-accessible functorial factorization that it suffices to verify that the full subcategory of $\CC^1$ comprised of trivial fibrations is a $\kappa$-accessibly embedded, $\kappa$-accessible subcategory. For this, consider the functor
\begin{equation*}
\Mor_{\CC,\Box}:\xymatrix@1@C=18pt@R=0pt{
\CC^1\ar[r]&\mathbf{Set}_{\XX}^{(I\cdot 1)}\\
[Y\to X]\,\ar@{|->}[r]&([\Mor_{\CC}(L,Y)\to\Mor_{\CC}(L,X)\times_{\Mor_{\CC}(K,X)}\Mor_{\CC}(K,Y)])_{[K\to L]\in I}.
}
\end{equation*}
Since the domains and codomains of $I$ are $\kappa$-presentable, this is a $\kappa$-accessible functor. The trivial fibrations are by definition the inverse image of the full subcategory of $I$-tuples of surjective morphisms under $\Mor_{\CC,\Box}$. 
\end{proof}
\end{prp}

\begin{cor} Any $\XX$-combinatorial model $\XX$-category satisfies the hypotheses of \ref{lem_smith-beke}.
\end{cor}

\begin{cor}\label{cor:tractislefttract} An $\XX$-combinatorial model $\XX$-category $\CC$ is $\XX$-tractable if and only if the $\XX$-small set $I$ of generating cofibrations can be chosen with cofibrant domains.
\begin{proof} Suppose $I$ is an $\XX$-small set of generating cofibrations with cofibrant domains, and suppose $J$ an $\XX$-small set of trivial cofibrations satisfying the conditions of \ref{lem:denseJ}. To give another such $\XX$-small set of trivial cofibrations with cofibrant domains, it suffices to show that any commutative square
\begin{equation*}
\xymatrix@C=18pt@R=18pt{
K\ar[d]\ar[r]&M\ar[d]\\
L\ar[r]&N
}
\end{equation*}
in which $[\fromto{K}{L}]\in I$ and $[\fromto{M}{N}]\in J$ can be factored as a a commutative diagram
\begin{equation*}
\xymatrix@C=18pt@R=18pt{
K\ar[d]\ar[r]&M'\ar[d]\ar[r]&M\ar[d]\\
L\ar[r]&N'\ar[r]&N,
}
\end{equation*}
in which $M'$ is cofibrant and $[\fromto{M'}{N'}]\in W\cap\cof I$. To construct this factorization, factor the morphism $\fromto{K}{M}$ as a cofibration $\fromto{K}{M'}$ followed by a weak equivalence $\fromto{M'}{M}$. Then factor the morphism $\fromto{L\sqcup^KM'}{N}$ as a cofibration $\fromto{L\sqcup^KM'}{N'}$ followed by a weak equivalence $\fromto{N'}{N}$. Then the composite $\fromto{M'}{N'}$ is a trivial cofibration providing the desired factorization.
\end{proof}
\end{cor}

\begin{nul} Suppose $\CC$ an $\XX$-combinatorial model $\XX$-category and $\DD$ a locally $\XX$-presentable category equipped with an adjunction
\begin{equation*}
F:\xymatrix@1@C=18pt{\CC\ar@<0.5ex>[r]&\DD\ar@<0.5ex>[l]}:U.
\end{equation*}
I now discuss circumstances under which the model structure on $\CC$ may be lifted to $\DD$.
\end{nul}

\begin{dfn}\label{dfn:projmodstr} 
\begin{enumerate}[(\ref{dfn:projmodstr}.1)]
\item A morphism $f:\fromto{X}{Y}$ of $\DD$ is said to be a \emph{projective weak equivalence} (respectively, a \emph{projective fibration}, a \emph{projective trivial fibration}) if $Uf:\fromto{UX}{UY}$ is a weak equivalence (resp., fibration, trivial fibration).
\item A morphism $f:\fromto{X}{Y}$ of $\DD$ is said to be a \emph{projective cofibration} if it satisfies the left lifting property with respect to any projective trivial fibration; $f$ is said to be a \emph{projective trivial cofibration} if it is, in addition, a projective weak equivalence.
\item If the projective weak equivalences, projective cofibrations, and projective fibrations define a model structure on $\DD$, then I call this model structure the \emph{projective model structure.}
\end{enumerate}
\end{dfn}

\begin{lem}\label{lem:Kan} Suppose that in $\DD$, transfinite compositions and pushouts of projective trivial cofibrations of $\DD$ are projective weak equivalences. Then the projective model structure on $\DD$ exists; it is $\XX$-combinatorial, and it is $\XX$-tractable if $\CC$ is. Furthermore the adjunction $(F,U)$ is a Quillen adjunction.
\begin{proof} The full accessible inverse image of an accessibly embedded accessible full subcategory is again an accessibly embedded accessible full subcategory; hence the projective weak equivalences are an accessibly embedded accessible subcategory of $\DD^1$. Choose now an $\XX$-small set $I$ of $\CC$ (respectively, $\CC_c$) of generating cofibrations.

One now applies the recognition lemma \ref{lem_smith-beke} to the set $W$ of projective weak equivalences and the $\XX$-small set $FI$. It is clear that $\inj I\subset W$, and by assumption it follows that $W\cap\cof I$ is closed under pushouts and transfinite compositions. One now verifies easily that the fibrations are the projective ones and that the adjunction $(F,U)$ is a Quillen adjunction.

Since $F$ is left Quillen, the set $FI$ has cofibrant domains if $I$ does.
\end{proof}
\end{lem}

\subsection*{Application I: Model structures on diagram categories} Suppose $\XX$ a universe, $K$ an $\XX$-small category, and $\CC$ an $\XX$-combinatorial (respectively, $\XX$-tractable) model $\XX$-category. The category $\CC(K)$ of $\CC$-valued presheaves on $K$ has two $\XX$-combinatorial (resp. $\XX$-tractable) model structures, to which we now turn.

\begin{dfn}\label{dfn:proj} A morphism $\fromto{X}{Y}$ of $\CC$-valued presheaves on $K$ is a \emph{projective weak equivalence} or \emph{projective fibration} if, for any object $k$ of $K$, the morphism $\fromto{X_k}{Y_k}$ is a weak equivalence or fibration of $\CC$.
\end{dfn}

\begin{thm} The category $\CC(K)$ of $\CC$-valued presheaves on $K$ has an $\XX$-combinatorial (resp., $\XX$-tractable) model structure --- the \emph{projective} model structure $\CC(K)_{\proj}$ ---, in which the weak equivalences and fibrations are the projective weak equivalences and fibrations.
\begin{proof} Consider the functor $e:\fromto{\Obj K}{K}$, which induces an adjunction
\begin{equation*}
e_!:\xymatrix@1@C=18pt{\CC(\Obj K)\ar@<0.5ex>[r]&\CC(K)\ar@<0.5ex>[l]}:e^{\star}.
\end{equation*}
The condition of \ref{lem:Kan} follows from the observation that $e^{\star}$ preserves all colimits.
\end{proof}
\end{thm}

\begin{dfn}\label{dfn:inj} A morphism $\fromto{X}{Y}$ of $\CC$-valued presheaves on $K$ is an \emph{injective weak equivalence} or \emph{injective cofibration} if, for any object $k$ of $K$, the morphism $\fromto{X_k}{Y_k}$ is a weak equivalence or cofibration of $\CC$.
\end{dfn}

\begin{thm} The category $\CC(K)$ of $\CC$-valued presheaves on $K$ has an $\XX$-combinatorial model structure --- the \emph{injective} model structure $\CC(K)_{\inj}$ ---, in which the weak equivalences and cofibrations are the injective weak equivalences and cofibrations.
\begin{proof} Suppose $\kappa$ an $\XX$-small regular cardinal such that $K$ is $\kappa$-small, $\CC$ is locally $\kappa$-presentable, and a set of generating cofibrations $I_{\CC}$ for $\CC$ can be chosen from $\CC_{\kappa}$ (resp., from $\CC_{\kappa}\cap\CC_{\lambda}$); without loss of generality, we may assume that $I_{\CC}$ is the $\XX$-small set of all cofibrations in $\CC_{\kappa}$ (resp., in $\CC_{\kappa}\cap\CC_{\lambda}$). Denote by $I_{\CC(K)}$ the set of injective cofibrations between $\kappa$-presentable objects of $\CC(K)$ (resp., between $\kappa$-presentable objects of $\CC(K)$ that are in addition objectwise cofibrant). This set contains a generating set of cofibrations for the projective model structure, so it follows that $\inj I_{\CC(K)}\subset W$.

The claim is now that any injective cofibration can be written as a retract of transfinite composition of pushouts of elements of $I_{\CC(K)}$. This point follows from a cardinality argument, which proceeds almost exactly as for $s\mathbf{Set}$-functors from an $s\mathbf{Set}$-category to a simplicial model category. For this cardinality argument I refer to \cite[A.3.3.12-14]{lurie_inftytopoi} of J. Lurie, whose proofs and exposition I am unable to improve upon.

Since colimits are formed objectwise, it follows that the injective trivial cofibrations are closed under pushouts and transfinite composition.
\end{proof}
\end{thm}

\begin{prp} The identity functor is a Quillen equivalence
\begin{equation*}
\xymatrix@1@C=18pt{\CC(K)_{\proj}\ar@<0.5ex>[r]&\CC(K)_{\inj}\ar@<0.5ex>[l]}.
\end{equation*}
\begin{proof} Projective cofibrations are in particular objectwise cofibrations, and weak equivalences are identical sets.
\end{proof}
\end{prp}

\begin{prp} If $\CC$ is left or right proper, then so are $\CC(K)_{\proj}$ and $\CC(K)_{\inj}$.
\begin{proof} Pullbacks and pushouts are defined objectwise; hence it suffices to note that in both model structures, weak equivalences are defined objectwise, and any cofibration or fibration is in particular an objectwise cofibration or fibration.
\end{proof}
\end{prp}

\begin{prp} A functor $f:\fromto{K}{L}$ induces Quillen adjunctions
\begin{equation*}
f_!:\xymatrix@1@C=18pt{\CC(K)_{\proj}\ar@<0.5ex>[r]&\CC(L)_{\proj}\ar@<0.5ex>[l]}:f^{\star}\textrm{\qquad and\qquad}f^{\star}:\xymatrix@1@C=18pt{\CC(L)_{\inj}\ar@<0.5ex>[r]&\CC(K)_{\inj}\ar@<0.5ex>[l]}:f_{\star},
\end{equation*}
which are of course equivalences of categories if $f$ is an equivalence of categories.
\begin{proof} Clearly $f^{\star}$ preserves objectwise weak equivalences, objectwise cofibrations, and objectwise fibrations.
\end{proof}
\end{prp}

\subsection*{Application II: Model structures on section categories} Left and right Quillen presheaves are diagrams of model categories. Here I give model structures on their categories of sections, analogous to the injective and projective model structures on diagram categories above.

\begin{nul} Suppose here $\XX$ a universe.
\end{nul}

\begin{dfn}\label{dfn:ltrtQpresheaves} Suppose $K$ an $\XX$-small category.
\begin{enumerate}[(\ref{dfn:ltrtQpresheaves}.1)]
\item A \emph{left} (respectively, \emph{right}) \emph{Quillen presheaf} on $K$ is a functor\footnote{In practice, of course, one is usually presented with a pseudofunctor, rather than a functor. Well-known rectification results allow one to replace such a pseudofunctor with a pseudoequivalent functor, and all the model structures can be lifted along this pseudoequivalence.} $\FF:\fromto{K^{\op}}{\mathbf{Cat}_{\YY}}$ for some universe $\YY$ with $\XX\in\YY$ such that for every $k\in\Obj K$, the category $\FF_k$ is a model $\XX$-category, and for every morphism $f:\fromto{\ell}{k}$ of $K$, the induced functor $f^{\star}:\fromto{\FF_k}{\FF_{\ell}}$ is left (resp., right) Quillen.
\item A left or right Quillen presheaf $\FF$ on $K$ is said to be \emph{$\XX$-combinatorial} (respectively, \emph{$\XX$-tractable}, \emph{left proper}, \emph{right proper}, ...) if for every $k\in\Obj K$, the model $\XX$-category $\FF_k$ is so.
\item A \emph{left} (respectively, \emph{right}) \emph{morphism} $\Theta:\fromto{\FF}{\GG}$ \emph{of left} (resp., \emph{right}) \emph{Quillen presheaves} is a pseudomorphism of functors $\fromto{K^{\op}}{\mathbf{Cat}_{\YY}}$ such that for any $k\in\Obj K$, the functor $\Theta_k:\fromto{\FF_k}{\GG_k}$ is left (resp., right) Quillen.
\item A \emph{left} (respectively, \emph{right}) \emph{section} $X$ of a left (resp., right) Quillen presheaf $\FF$ is a tuple $(X,\phi)=((X_k)_{k\in\Obj K},(\phi_f)_{f\in\Obj(K^1)})$ comprised of an object $X=(X_k)_{k\in\Obj K}$ of $\Prod_{k\in\Obj K}\FF_k$ and a morphism $\phi_f:\fromto{f^{\star}X_{\ell}}{X_k}$ (resp., $\phi_f:\fromto{X_k}{f^{\star}X_{\ell}}$), one for each morphism $[f:\fromto{\ell}{k}]\in K$, such that for any composable pair
\begin{equation*}
[\xymatrix@1@C=18pt{m\ar[r]^g&\ell\ar[r]^f&k}]\in K,
\end{equation*}
one has the identity
\begin{equation*}
\phi_g\circ(g^{\star}\phi_f)=\phi_{f\circ g}:\fromto{(f\circ g)^{\star}X_k}{X_m}
\end{equation*}
(resp.,
\begin{equation*}
(g^{\star}\phi_f)\circ\phi_g=\phi_{f\circ g}:\fromto{X_m}{(f\circ g)^{\star}X_k}\textrm{\quad).}
\end{equation*}
\item A \emph{morphism of left} (respectively, \emph{right}) \emph{sections} $r:\fromto{(X,\phi)}{(Y,\psi)}$ is a morphism $r:\fromto{X}{Y}$ of $\Prod_{k\in\Obj K}\FF_k$ such that the diagram
\begin{equation*}
\xymatrix@C=24pt@R=24pt{
f^{\star}X_k\ar[d]_{\phi_f}\ar[r]^{f^{\star}r_k}&f^{\star}Y_k\ar[d]^{\psi_f}\\
X_{\ell}\ar[r]_{r_{\ell}}&Y_{\ell}
}
\end{equation*}
in $\FF_{\ell}$ (resp., the diagram
\begin{equation*}
\xymatrix@C=24pt@R=24pt{
X_k\ar[d]_{\phi_f}\ar[r]^{r_k}&Y_k\ar[d]^{\psi_f}\\
f^{\star}X_{\ell}\ar[r]_{f^{\star}r_{\ell}}&f^{\star}Y_{\ell}
}
\end{equation*}
in $\FF_k$) commutes for any morphism $f:\fromto{\ell}{k}$ of $K$. These morphisms clearly compose, to give a category $\Sect^L\FF$ (resp., $\Sect^R\FF$) of left (resp., right) sections of $\FF$.
\end{enumerate}
\end{dfn}

\begin{lem} A left (respectively, right) morphism $\Theta:\fromto{\FF}{\GG}$ of left (resp., right) Quillen presheaves on an $\XX$-small category $K$ induces a left adjoint
\begin{equation*}
\Theta_!:\fromto{\Sect^L\FF}{\Sect^L\GG}
\end{equation*}
(resp., a right adjoint
\begin{equation*}
\Theta_{\star}:\fromto{\Sect^L\FF}{\Sect^L\GG}\textrm{\quad).}
\end{equation*}
\begin{proof} If $\Theta$ is a left morphism of left Quillen presheaves, then define $\Theta_!$ by the formula
\begin{equation*}
\Theta_!(X,\phi):=((\Theta_kX_k)_{k\in\Obj K},((\Theta\phi_f)\circ\theta_f)_{f\in\Obj(K^1)})
\end{equation*}
for any left section $(X,\phi)=((X_k)_{k\in\Obj K},(\phi_f)_{f\in\Obj(K^1)})$, in which the morphism
\begin{equation*}
\theta_f:\fromto{f^{\star}\Theta_kX_k}{\Theta_{\ell}f^{\star}X_k}
\end{equation*}
is the structural isomorphism of the pseudomorphism $\Theta$.

Its right adjoint
\begin{equation*}
\Theta^{\star}:\fromto{\Sect^L\GG}{\Sect^L\FF}
\end{equation*}
is defined by the formula
\begin{equation*}
\Theta^{\star}(Y,\psi):=((H_kY_k)_{k\in\Obj K},(\eta_f)_{f\in\Obj(K^1)})
\end{equation*}
for any left section $(Y,\psi)=((Y_k)_{k\in\Obj K},(\psi_f)_{f\in\Obj(K^1)})$, in which the morphism
\begin{equation*}
\eta_f:\fromto{f^{\star}H_kY_k}{H_{\ell}f^{\star}Y_k}
\end{equation*}
is the morphism adjoint to the composite
\begin{equation*}
\xymatrix@1@C=24pt{\Theta_{\ell}f^{\star}H_kY_k\ar[r]^-{\theta_f}&f^{\star}\Theta_kH_kY_k\ar[r]^-{f^{\star}c}&f^{\star}H_k},
\end{equation*}
where $c:\fromto{\Theta_kH_kY_k}{Y_k}$ is the counit of the adjunction $(\Theta_k,H_k)$.

The corresponding statement for right morphisms follows by duality.
\end{proof}
\end{lem}

\begin{nul} The previous lemma suggests that the most natural model structure on the category of left (respectively, right) sections of a left (resp., right) Quillen presheaf $\FF$ is an injective (resp., projective) one, in which the weak equivalences and cofibrations (resp., fibrations) are defined objectwise. This idea is borne out by the observation that these model categories can be thought of as good models for the $(\infty,1)$-categorical lax limit (resp., $(\infty,1)$-categorical colax limit) of $\FF$.
\end{nul}

\begin{lem} If $\FF$ is a left (respectively, right) $\XX$-combinatorial Quillen presheaf on an $\XX$-small category $K$, then the category $\Sect^L\FF$ (resp., $\Sect^R\FF$) is locally $\XX$-presentable.
\begin{proof} It is a simple matter to verify that $\Sect^L\FF$ and $\Sect^R\FF$ are complete and cocomplete. The category of left sections is the lax limit of $\FF$, and the category of right sections is a colax limit of $\FF$; so the result follows from the fact that the $2$-category of $\XX$-accessible categories is closed under arbitrary $\XX$-small weighted bilimits in which all functors are accessible \cite[Theorem 5.1.6]{MR1031717}.
\end{proof}
\end{lem}

\begin{lem} Suppose $a:\fromto{L}{K}$ a functor of $\XX$-small categories. If $\FF$ is a left (respectively, right) $\XX$-combinatorial Quillen presheaf on an $\XX$-small category $L$, then there is a string $(a_!,a^{\star},a_{\star})$ of adjoints
\begin{equation*}
a_!,a_{\star}:\xymatrix@1@C=18pt{\Sect^L(\FF\circ a)\ar@<1ex>[r]\ar@<-1ex>[r]&\Sect^L\FF\ar[l]}:a^{\star}\textrm{\quad (resp.,\quad}a_!,a_{\star}:\xymatrix@1@C=18pt{\Sect^R(\FF\circ a)\ar@<1ex>[r]\ar@<-1ex>[r]&\Sect^R\FF\ar[l]}:a^{\star}\textrm{\quad).}
\end{equation*}
\begin{proof} If $\FF$ is a left Quillen presheaf, then the functor
\begin{equation*}
a^{\star}:\fromto{\Sect^L(\FF\circ a)}{\Sect^L\FF}
\end{equation*}
is simply given by the formula
\begin{equation*}
a^{\star}(X,\phi):=((X_{a(k)})_{k\in\Obj K},(\phi_{a(f)})_{f\in\Obj(K^1)}).
\end{equation*}
Since this functor commutes with all limits and colimits, the existence of its left and right adjoints follows from the usual adjoint functor theorems.
\end{proof}
\end{lem}

\begin{dfn} Suppose $K$ an $\XX$-small category, $\FF$ a right Quillen presheaf on $K$. A morphism $\fromto{X}{Y}$ of right sections of $\FF$ is a \emph{projective weak equivalence} or \emph{projective fibration} if, for any object $k$ of $K$, the morphism $\fromto{X_k}{Y_k}$ is a weak equivalence or fibration of $\FF_k$.
\end{dfn}

\begin{thm} The category $\Sect^R\FF$ of right sections of an $\XX$-combinatorial (respectively, $\XX$-tractable) right Quillen presheaf $\FF$ on an $\XX$-small category $K$ has an $\XX$-combinatorial (resp., $\XX$-tractable) model structure --- the \emph{projective} model structure $\Sect^R_{\proj}\FF$ ---, in which the weak equivalences and fibrations are the projective weak equivalences and fibrations.
\begin{proof} Consider the functor $e:\fromto{\Obj K}{K}$, which induces an adjunction
\begin{equation*}
e_!:\xymatrix@1@C=18pt{\Prod_{k\in\Obj K}\FF_k\ar@<0.5ex>[r]&\Sect^R\FF\ar@<0.5ex>[l]}:e^{\star}.
\end{equation*}
The condition of \ref{lem:Kan} follows from the observation that $e^{\star}$ preserves all colimits.
\end{proof}
\end{thm}

\begin{dfn} Suppose $K$ an $\XX$-small category, $\FF$ a left Quillen presheaf on $K$. A morphism $\fromto{X}{Y}$ of left sections of $\FF$ is an \emph{injective weak equivalence} or \emph{injective cofibration} if, for any object $k$ of $K$, the morphism $\fromto{X_k}{Y_k}$ is a weak equivalence or cofibration of $\FF_k$.
\end{dfn}

\begin{thm} The category $\Sect^L\FF$ of left sections of an $\XX$-combinatorial (respectively, $\XX$-tractable) left Quillen presheaf $\FF$ on an $\XX$-small category $K$ has an $\XX$-combinatorial (resp., $\XX$-tractable) model structure --- the \emph{injective} model structure $\Sect^L_{\inj}\FF$ ---, in which the weak equivalences and cofibrations are the injective weak equivalences and cofibrations.
\begin{proof} Suppose $\kappa$ an $\XX$-small regular cardinal such that $K$ is $\kappa$-small, each $\FF_k$ is locally $\kappa$-presentable, and a set of generating cofibrations $I_{\FF_k}$ for each $\FF_k$ can be chosen from $\FF_{k,\kappa}$ (resp., from $\FF_{k,\kappa}\cap\FF_{k,c}$); without loss of generality, we may assume that $I_{\FF_k}$ is the $\XX$-small set of all cofibrations in $\FF_{k,\kappa}$ (resp., in $\FF_{k,\kappa}\cap\FF_{k,c}$). Denote by $I_{\Sect^L\FF}$ the set of injective cofibrations between $\kappa$-presentable objects of $\Sect^L\FF$ (resp., between $\kappa$-presentable objects of $\Sect^L\FF$ that are in addition objectwise cofibrant). This set contains a generating set of cofibrations for the projective model structure, so it follows that $\inj I_{\Sect^L\FF}\subset W$.

The argument given for the existence of the injective model structure on presheaf categories applies almost verbatim here to demonstrate that any injective cofibration can be written as a retract of transfinite composition of pushouts of elements of $I_{\Sect^L\FF}$.

Since colimits are formed objectwise, it follows that the injective trivial cofibrations are closed under pushouts and transfinite composition.
\end{proof}
\end{thm}

\begin{prp} Suppose $K$ an $\XX$-small category, $\FF$ an $\XX$-combinatorial left (respectively, right) Quillen presheaf on $K$. If each $\FF_k$ is left or right proper, then so is $\Sect^L_{\inj}\FF$ (resp., $\Sect^R_{\proj}\FF$).
\begin{proof} Pullbacks and pushouts are defined objectwise; hence it suffices to note that in both model structures, weak equivalences are defined objectwise, and any cofibration or fibration is in particular an objectwise cofibration or fibration.
\end{proof}
\end{prp}

\begin{prp} A left (resp., right) morphism $\Theta:\fromto{\FF}{\GG}$ of $\XX$-combinatorial left (resp., right) Quillen presheaves on an $\XX$-small category $K$ induces a Quillen adjunction
\begin{equation*}
\Theta_!:\xymatrix@1@C=18pt{\Sect^L_{\inj}\FF\ar@<0.5ex>[r]&\Sect^L_{\inj}\GG\ar@<0.5ex>[l]}:\Theta^{\star}\textrm{\quad (resp.,\quad}\Theta^{\star}:\xymatrix@1@C=18pt{\Sect^R_{\proj}\GG\ar@<0.5ex>[r]&\Sect^R_{\proj}\FF\ar@<0.5ex>[l]}:\Theta_{\star}\textrm{\quad),}
\end{equation*}
which is a Quillen equivalence if each $\Theta_k$ is a Quillen equivalence.
\begin{proof} Clearly $\Theta_!$ (resp., $\Theta_{\star}$) preserves objectwise weak equivalences and objectwise cofibrations (resp., fibrations).
\end{proof}
\end{prp}

\section{Left Bousfield localization} 

\subsection*{Definition and existence of left Bousfield localizations} Here I review some highlights from the general theory of left Bousfield localization. Since the published references do not include a proof of the existence theorem of J. Smith, I include it solely for convenience of reference, with the caveat that the result should in no way be construed as mine.

\begin{nul} Suppose $\XX$ a universe, $\MM$ a model $\XX$-category.
\end{nul}

\begin{nul} Suppose $X$ and $Y$ objects of $\MM$. For any odd natural number $n$, I define a category $w\MOR_{\MM}^n(X,Y)$. The objects of $w\MOR_{\MM}^n(X,Y)$ are strings of morphisms
\begin{equation*}
\xymatrix@C=18pt{X=X_0&\ar[l]X_1\ar[r]&X_2&\ar[l]\dots\ar[r]&X_{n-1}&\ar[l]X_n=Y}
\end{equation*}
such that each morphism $\fromto{X_{2i}}{X_{2i+1}}$ is contained in $w\MM$. Morphisms between two such sequences are simply commutative diagrams of the form
\begin{equation*}
\xymatrix@C=18pt@R=9pt{
&\ar[dl]X_1\ar[dd]\ar[r]&X_2\ar[dd]&\ar[l]\dots\ar[r]&X_{n-1}\ar[dd]&\\
X&&&&&Y\ar[dl]\ar[ul]\\
&\ar[ul]X'_1\ar[r]&X'_2&\ar[l]\dots\ar[r]&X'_{n-1}&,
}
\end{equation*}
wherein the vertical maps are in $w\MM$. Recall that the \emph{hammock localization} of $\MM$ is the $s\mathbf{Set}$-category $L^H\MM$ whose objects are exactly those of $\MM$, with
\begin{equation*}
\MOR_{L^H\MM}(X,Y)=\colim_n\nu_{\bullet}(w\MOR_{\MM}^n(X,Y))
\end{equation*}
for any objects $X$ and $Y$.

The standard references on the hammock localization are the triple of papers \cite{MR81h:55019}, \cite{MR81m:55018}, and \cite{MR81h:55018} of W. G. Dwyer and D. Kan. A more modern treatment can be found in \cite{MR2102294}.
\end{nul}

\begin{sch}[Dwyer--Kan] Suppose $Q:\fromto{\MM}{\MM_c}$ a cofibrant replacement functor, $R:\fromto{\MM}{\MM_f}$ a fibrant replacement functor, $\Gamma^{\bullet}:\fromto{\MM}{(c\MM)_c}$ a cosimplicial resolution functor, and $\Lambda_{\bullet}:\fromto{\MM}{(s\MM)_f}$ a simplicial resolution functor; then there are natural weak equivalences of the simplicial sets
\begin{equation*}
\xymatrix@C=-36pt@R=18pt{
\Mor_{\MM}(\Gamma^{\bullet}X,RY)\ar[dr]&&\Mor_{\MM}(QX,\Lambda_{\bullet}Y)\ar[dl]\\
&\diag\Mor_{\MM}(\Gamma^{\bullet}X,\Lambda_{\bullet}Y)&\\
&\ar[u]\hocolim_{(p,q)\in\Delta^{\op}\times\Delta^{\op}}\Mor_{\MM}(\Gamma^pX,\Lambda_qY)\ar[d]&\\
&\nu_{\bullet}w\Mor_{\MM}^3(X,Y)\ar[d]&\\
&\Mor_{L^H\MM}(X,Y).&
}
\end{equation*}
\end{sch}

\begin{ntn}\label{ntn_modcat} Write $\RMor_{\MM}$ for the simplicial mapping space functor
\begin{equation*}
\xymatrix@C=10pt@R=2pt{\Ho\MM^{\op}\times\Ho\MM\ar[r]&\Ho s\mathrm{Set}_{\XX}\\
(X,Y)\,\ar@{|->}[r]&\MOR_{L\MM}(X,Y);}
\end{equation*}
the subscript may omitted when the context is sufficiently clear.
\end{ntn}

\begin{lem}[\protect{\cite[17.7.7]{MR2003j:18018}}]\label{lem_RMordetectsws} The following are equivalent for a morphism $\fromto{A}{B}$ of $\MM$.
\begin{enumerate}[(\ref{lem_RMordetectsws}.1)]
\item The morphism $\fromto{A}{B}$ is a weak equivalence.
\item For any fibrant object $Z$ of $\MM$, the induced morphism
\begin{equation*}
\fromto{\RMor_{\MM}(B,Z)}{\RMor_{\MM}(A,Z)}
\end{equation*}
is an isomorphism of $\Ho s\mathrm{Set}_{\XX}$.
\item For any cofibrant object $X$ of $\MM$, the induced morphism
\begin{equation*}
\fromto{\RMor_{\MM}(X,A)}{\RMor_{\MM}(X,B)}
\end{equation*}
is an isomorphism of $\Ho s\mathrm{Set}_{\XX}$.
\end{enumerate}
\end{lem}

\begin{dfn}\label{dfn_hoclass} Suppose $H$ a set of homotopy classes of morphisms of $\MM$.
\begin{enumerate}[(\ref{dfn_hoclass}.1)]
\item A \emph{left Bousfield localization of $\MM$ with respect to $H$} is a model $\XX$-category $L_H\MM$, equipped with a left Quillen functor $\fromto{\MM}{L_H\MM}$ that is initial among left Quillen functors $F:\fromto{\MM}{\NN}$ to model $\XX$-categories $\NN$ with the property that for any $f$ representing a class in $H$, $\LL F(f)$ is an isomorphism of $\Ho\NN$.
\item An object $Z$ of $\MM$ is \emph{$H$-local} if for any morphism $\fromto{A}{B}$ representing an element of $H$, the morphism
\begin{equation*}
\fromto{\RMor_{\MM}(B,Z)}{\RMor_{\MM}(A,Z)}
\end{equation*}
is an isomorphism of $\Ho s\mathbf{Set}_{\XX}$.
\item A morphism $\fromto{A}{B}$ of $\MM$ is an \emph{$H$-local equivalence} if for any $H$-local object $Z$, the morphism
\begin{equation*}
\fromto{\RMor_{\MM}(B,Z)}{\RMor_{\MM}(A,Z)}
\end{equation*}
is an isomorphism of $\Ho s\mathbf{Set}_{\XX}$.
\end{enumerate}
\end{dfn}

\begin{lem} When it exists, the left Bousfield localization of $\MM$ with respect to $H$ is unique up to a unique isomorphism of model $\XX$-categories.
\begin{proof} Initial objects are essentially unique.
\end{proof}
\end{lem}

\begin{nul} Left Bousfield localizations of left proper, $\XX$-combinatorial model $\XX$-categories with respect to $\XX$-small sets of homotopy classes of morphisms are guaranteed to exist, as I shall now demonstrate. For the remainder of this section, suppose $H$ a set of homotopy classes of morphisms of $\MM$. The uniqueness, characterization, and existence of left Bousfield localizations are the central objectives of the next few results. Uniqueness is a simple matter, and it is a familiar fact that if a model structure exists on $\MM$ with the same cofibrations whose weak equivalences are the $H$-local weak equivalences, then this is the left Bousfield localization. The central point is thus to determine the existence of such a model structure. Smith's existence theorem \ref{thm_lbexist} hinges on the recognition principle \ref{lem_smith-beke} and the following pair of technical lemmata.
\end{nul}

\begin{lem} If $\MM$ is $\XX$-combinatorial, and $H$ is $\XX$-small, then the set of $\HH$-local objects of $\MM$ comprise an accessibly embedded, accessible subcategory of $\MM^1$.
\begin{proof}  Choose an accessible fibrant replacement functor $R$ for $\MM$, a functorial cosimplicial resolution functor $\Gamma^{\bullet}:\fromto{\MM}{(c\MM)_c}$, and an $\XX$-small set $S$ of representatives for all and only the homotopy classes of $H$. Then the functor
\begin{equation*}
\xymatrix@1@C=18pt@R=0pt{\MM\ar[r]&(s\mathbf{Set}_{\XX})^1\\
Z\,\ar@{|->}[r]&\Coprod_{f\in S}f_{\Gamma^{\bullet},R}^{\star}(Z)}
\end{equation*}
is accessible, where $f_{\Gamma^{\bullet},R}^{\star}(Z):\fromto{\Mor_{\MM}(\Gamma^{\bullet} Y,RZ)}{\Mor_{\MM}(\Gamma^{\bullet}X,RZ)}$ is the morphism of simplicial sets induced by $f:\fromto{X}{Y}$. Since the full subcategory of $(s\mathbf{Set}_{\XX})^1$ comprised of weak equivalences is accessibly embedded and accessible, the full subcategory of $H$-local objects is also accessibly embedded and accessible.
\end{proof}
\end{lem}

\begin{lem}\label{lem_Hlocequivaccessible} If $\MM$ is $\XX$-combinatorial, and $H$ is $\XX$-small, then the set of $H$-local equivalences of $\MM$ comprise an accessibly embedded, accessible subcategory of $\MM^1$.
\begin{proof} This follows from the previous lemma and the fact that for sufficiently large regular $\XX$-small cardinals $\kappa$, the $H$-local equivalences of $\MM$ are closed under $\kappa$-filtered colimits. To show the latter point, choose $\kappa$ so that $\kappa$-filtered colimits are homotopy colimits. Then for any $H$-local object $Z$, a colimit $\fromto{\colim A}{\colim B}$ of a $\kappa$-filtered diagram of $H$-local equivalences is a weak equivalence because the morphism
\begin{equation*}
\fromto{\RMor_{\MM}(\colim B,Z)}{\RMor_{\MM}(\colim A,Z)}
\end{equation*}
is a homotopy limit of weak equivalences in $s\mathrm{Set}_{\XX}$, hence a weak equivalence.
\end{proof}
\end{lem}

\begin{thm}[Smith, \protect{\cite{combinatorial}}]\label{thm_lbexist} If $\MM$ is $\XX$-combinatorial, and $H$ is an $\XX$-small set of homotopy classes of morphisms of $\MM$, the left Bousfield localization $L_H\MM$ of $\MM$ along any set representing $H$ exists and satisfies the following conditions.
\begin{enumerate}[(\ref{thm_lbexist}.1)]
\item The model category $L_H\MM$ is left proper and $\XX$-combinatorial.
\item As a category, $L_H\MM$ is simply $\MM$.
\item\label{item_cofsLM} The cofibrations of $L_H\MM$ are exactly those of $\MM$.
\item The fibrant objects of $L_H\MM$ are the fibrant $H$-local objects $Z$ of $\MM$.
\item\label{item_wesLM} The weak equivalences of $L_H\MM$ are the $H$-local equivalences.
\end{enumerate}
\begin{proof} The aim is to guarantee that a cofibrantly generated model structure on $\MM$ exists sstisfying conditions (\ref{thm_lbexist}.\ref{item_cofsLM})--(\ref{thm_lbexist}.\ref{item_wesLM}) using \ref{lem_smith-beke}. The combinatoriality is then automatic, and the universal property and the left properness are then verified in \cite[Theorem 3.3.19 and Proposition 3.4.4]{MR2003j:18018}.

Fix an $\XX$-small set $I_{\MM}$ of generating cofibrations of $\MM$, and let $wL_H\MM$ denote the set of the weak equivalences described in (\ref{thm_lbexist}.\ref{item_wesLM}). By \ref{lem_Hlocequivaccessible}, we can now apply \ref{lem_smith-beke}: observe that since $I_{\MM}$-injectives are trivial fibrations of $\MM$, they are in particular weak equivalences of $\MM$, and hence are among the elements of $wL_H\MM$.

It thus remains only to show that pushouts and transfinite compositions of morphisms of $\cof\MM\cap wL_H\MM$ are $H$-local weak equivalences. Suppose first that $\fromto{K}{L}$ a cofibration in $wL_H\MM$, and suppose
\begin{equation*}
\xymatrix@C=18pt@R=18pt{
K\ar[d]\ar[r]&L\ar[d]\\
K'\ar[r]&L'
}
\end{equation*}
a pushout diagram in $\MM$. Note that by the left properness of $\MM$, this pushout is in fact a homotopy pushout; thus the statement that $\fromto{K'}{L'}$ is an element of $wL_H\MM$ is equivalent to the assertion that, for any $H$-local object $Z$, the diagram
\begin{equation*}
\xymatrix@C=18pt@R=18pt{
\RMor_{\MM}(K',Z)\ar[d]\ar[r]&\RMor_{\MM}(K',Z)\ar[d]\\
\RMor_{\MM}(L,Z)\ar[r]&\RMor_{\MM}(K,Z)
}
\end{equation*}
is a homotopy pullback diagram in $s\mathbf{Set}$, and this follows immediately from the fact that $\fromto{\RMor_{\MM}(L,Z)}{\RMor_{\MM}(K,Z)}$ is a weak equivalence. Since $\kappa$-filtered colimits are homotopy colimits for $\kappa$ sufficiently large, it follows that a transfinite composition of elements of $\cof\MM\cap wL_H\MM$ is an morphism of $wL_H\MM$.
\end{proof}
\end{thm}

\begin{dfn} If $\MM$ is left proper and $\XX$-combinatorial, and if $H$ is an $\XX$-small set of homotopy classes of morphisms of $\MM$, then an object $X$ of the left Bousfield localization $L_H\MM$ is \emph{quasifibrant} if some fibrant replacement $R_{\MM}X$ of $X$ in $\MM$ is fibrant in $L_H\MM$.
\end{dfn}

\begin{lem} If $\MM$ is left proper and $\XX$-combinatorial, and $H$ is an $\XX$-small set of homotopy classes of morphisms of $\MM$, then an object $X$ of the left Bousfield localization $L_H\MM$ is quasifibrant if and only if it is $H$-local.
\begin{proof} $H$-locality is closed under weak equivalences in $\MM$; hence if $X$ is quasifibrant it is surely $H$-local, and the fibrant replacement in $\MM$ of an $H$-local object is $H$-local.
\end{proof}
\end{lem}

\begin{nul} As a rule, one has essentially no control in a left Bousfield localization over the generating trivial cofibrations. The following proposition (originally --- with a different proof --- due to M. Hovey) is one of the very few results on the trivial cofibrations of left Bousfield localizations; it is critical for the forthcoming existence theorem \ref{thm_elbexist} for \emph{enriched} left Bousfield localizations.
\end{nul}

\begin{prp}[Hovey, \protect{\cite[Proposition 4.3]{MR2066503}}]\label{prp_HoveyJLHMJ} Suppose that $\MM$ is left proper and $\XX$-combinatorial, and suppose that $H$ is $\XX$-small. Then the left Bousfield localization $L_H\MM$ is $\XX$-tractable if $\MM$ is.
\begin{proof} Immediate from \ref{cor:tractislefttract}.
\end{proof}
\end{prp}

\subsection*{The failure of right properness} Left Bousfield localizations inherit left properness, but in general they destroy right properness. This is because there is very little control over the fibrations.

Nevertheless, there are often full subcategories that are in a sense right proper, and this form of right properness is inherited by the quasifibrant objects contained in these subcategories in the left Bousfield localization. In this case, there exist functorial factorizations of morphisms of these quasifibrant objects through quasifibrant objects, so homotopy pullbacks of these quasifibrant objects can be computed effectively. This also provides a nice recognition principle for fibrations of $L_H\MM$ with a quasifibrant codomain that lies in such a subcategory.

One can think of this subsection as an enlargement of Reedy's observation that homotopy pullbacks of fibrant objects can be computed by replacing on only one side, or, alternatively, one can think of this subsection as a collection of techniques for coping with the reality that many important combinatorial model categories are simply not right proper.

There are, of course, dual conditions and results to many of those of this section, but they are not needed here, essentially because left properness is a relatively common condition in practice.

\begin{nul} Suppose $\XX$ a universe, $\MM$ a model $\XX$-category.
\end{nul}

\begin{dfn}\label{dfn_robust}\begin{enumerate}[(\ref{dfn_robust}.1)]
\item  If $E$ is any full subcategory of $\MM$, an \emph{$E$-placement functor} is a pair $(r_E,\epsilon_E)$ consisting of a functor $r_E:\fromto{\MM}{E}$ along with an objectwise weak equivalence $\epsilon_E$ from the identity functor to the composite $\iota_E\circ r_E$, where $\iota_E:\fromto{E}{\MM}$ denotes the inclusion.
\item A full subcategory $E$ of $\MM$ is said to be \emph{stable under trivial fibrations} if for any object $X$ of $E$ and any trivial fibration $\fromto{Y}{X}$, the object $Y$ is an object of $E$ as well.
\item A pullback diagram
\begin{equation*}
\xymatrix@C=18pt@R=18pt{Y'\ar[r]\ar[d]&Y\ar[d]\\
X'\ar[r]&X}
\end{equation*}
in $\MM$ in which $\fromto{Y}{X}$ is a fibration is called an \emph{admissible pullback diagram}, and the morphism $\fromto{Y'}{Y}$ is called the \emph{admissible pullback} of $\fromto{X'}{X}$ along $\fromto{Y}{X}$.
\item A full subcategory $E$ of $\MM$ is said to be \emph{admissibly left exact} if $E$ contains any admissible pullback in $E$ --- i.e., if $E$ contains the pullback of any morphism of $E$ along any fibration of $E$; more generally, $E$ is said to be \emph{partially left exact} if $E$ contains the pullback of any weak equivalence of $E$ along any fibration of $E$.
\item A partially left exact full subcategory $E$ is said to be \emph{right proper} if in $E$ admissible pullbacks of weak equivalences are weak equivalences.
\item Suppose $E$ a partially left exact full subcategory; then a pair $(F,\eta)$ consisting of a functor $F:\fromto{E}{\MM}$ with an objectwise weak equivalence $\eta$ from the forgetful functor to $F$ is said to be \emph{exceptional on $E$} if $F$ preserves admissible pullback diagrams.
\end{enumerate}
\end{dfn}

\begin{lem} Suppose $E$ and $E'$ partially (respectively, admissibly) left exact full subcategories of $\MM$ and $(F,\eta)$ an exceptional pair on $E$. Then the full subcategory $F^{-1}(E')$ of $E$ consisting of those objects $X$ of $E$ such that $FX$ is an object of $E'$ is partially (resp., admissibly) left exact. Moreover, if $E'$ is right proper, then so is $F^{-1}(E')$.
\begin{proof} Refer to the diagram
\begin{equation*}
\xymatrix@C=9pt@R=9pt{
F(Y')\ar[rrrr]\ar[dddd]&&&&FY\ar[dddd]\\
&Y'\ar[rr]\ar[dd]\ar[ul]&&Y\ar[dd]\ar[ur]&\\
&&&&\\
&X'\ar[rr]\ar[dl]&&X\ar[dr]&\\
F(X')\ar[rrrr]&&&&FX.
}
\end{equation*}
If the interior square is an admissible pullback diagram in $E$, then the outer square is so in $\MM$. If $\fromto{X'}{X}$ is a weak equivalence, then so is $\fromto{F(X')}{FX}$, and if $\fromto{F(Y')}{FY}$ is a weak equivalence, then so is $\fromto{Y'}{Y}$.
\end{proof}
\end{lem}

\begin{lem}\label{lem_ReedyMfrightproper} The full subcategory $\MM_f$ of fibrant objects is admissibly left exact and right proper.
\begin{proof} This follows from Reedy's argument, \cite[Theorem B]{reedy} or \cite[Proposition 13.1.2]{MR2003j:18018}.
\end{proof}
\end{lem}

\begin{cor} A partially left exact full subcategory $E$ of $\MM$ is right proper if there exists an exceptional fibrant replacement functor on $E$.
\end{cor}

\begin{cor} The model category $s\mathbf{Set}$ of simplicial sets is right proper.
\begin{proof} Kan's $\mathrm{Ex}^{\infty}$ is an exceptional fibrant replacement functor.
\end{proof}
\end{cor}

\begin{lem}\label{lem_Ertpropfunctfact} Suppose $E$ a right proper, partially left exact full subcategory of $\MM$ with $\MM_f\subset E$. Then there exists a functorial factorization of any morphism in $E$ into a weak equivalence of $E$ followed by a fibration of $E$. If, in addition, $E$ is stable under trivial fibrations, then there is a functorial factorization of every morphism of $E$ into a trivial cofibration of $E$ followed by a fibration of $E$.
\begin{proof} Choose a functorial factorization of every morphism into a trivial cofibration followed by a fibration; this gives in particular a functorial fibrant replacement $r$. Suppose $\fromto{X}{Y}$ a morphism of $E$; applying the chosen functorial factorization to the vertical morphism on the right in the diagram
\begin{equation*}
\xymatrix@C=18pt@R=18pt{
X\ar[r]\ar[d]&rX\ar[d]\\
Y\ar[r]&rY,
}
\end{equation*}
yields $\xymatrix@1@C=10pt{rX\ar[r]&Z\ar[r]&rY}$. Pulling back the resulting fibration $\fromto{Z}{rY}$ along $\fromto{Y}{rY}$ produces a diagram
\begin{equation*}
\xymatrix@C=18pt@R=18pt{
X\ar[r]\ar[d]&rX\ar[d]\\
Y\times_{rY}Z\ar[r]\ar[d]&Z\ar[d]\\
Y\ar[r]&rY,
}
\end{equation*}
in which $\fromto{Y\times_{rY}Z}{Y}$ is a fibration of $E$, and $\fromto{Y\times_{rY}Z}{Z}$ --- and therefore also $\fromto{X}{Y\times_{rY}Z}$ --- is a weak equivalence of $E$. If $E$ is stable under trivial fibrations, then applying the chosen factorization to the morphism $\fromto{X}{Y\times_{rY}Z}$ provides the desired factorization of the second half of the statement.
\end{proof}
\end{lem}

\begin{cor}\label{cor_Eadmisshtypull} If $E$ is a right proper, admissibly left exact full subcategory of $\MM$ with $\MM_f\subset E$, then admissible pullbacks in $E$ are homotopy pullbacks.
\end{cor}

\begin{nul} Suppose for the remainder of this section that the model category $\MM$ is left proper and $\XX$-combinatorial, that $H$ is $\XX$-small, and that $E$ is a right proper, admissibly left exact full subcategory of $\MM$ with $\MM_f\subset E$. Write $\mathrm{loc}_E(H)$ for the full subcategory of $E$ consisting of $H$-local objects, viewed as a full subcategory of the Bousfield localization $L_H\MM$; observe that $(L_H\MM)_f\subset\mathrm{loc}_E(H)$. The objective is to show that $\mathrm{loc}_E(H)$ is right proper and admissibly left exact, whence the effective computability of homotopy pullbacks of $H$-local objects of $E$.
\end{nul}

\begin{lem} The subcategory $\mathrm{loc}_E(H)$ is a right proper partially left exact full subcategory of $L_H\MM$, and if $E$ is stable under trivial fibrations, then so is $\mathrm{loc}_E(H)$.
\begin{proof} Weak equivalences in $L_H\MM$ between $H$-local objects are weak equivalences of $\MM$, and the trivial fibrations of $L_H\MM$ are exactly those of $\MM$.
\end{proof}
\end{lem}

\begin{lem} If $E$ is stable under trivial fibrations, then $\mathrm{loc}_E(H)$ admits a functorial factorization, within $\mathrm{loc}_E(H)$, of every morphism into a trivial cofibration of $\MM$ followed by a fibration of $L_H\MM$.
\begin{proof} Applying \ref{lem_Ertpropfunctfact}, there is a functorial factorization of every morphism into a trivial cofibration of $\mathrm{loc}_E(H)$ followed by a fibration of $\mathrm{loc}_E(H)$. But a trivial cofibration in $L_H\MM$ between $H$-local objects is a trivial cofibration in $\MM$ as well.
\end{proof}
\end{lem}

\begin{cor} If $E$ is stable under trivial fibrations, then a morphism of $\mathrm{loc}_E(H)$ is a fibration in $L_H\MM$ if and only if it is a fibration of $\MM$.
\begin{proof} One implication is obvious; the other is a consequence of the retract argument.
\end{proof}
\end{cor}

\begin{prp}\label{prp_locErpropadmisslex} The subcategory $\mathrm{loc}_E(H)$ is a right proper, admissibly left exact full subcategory of $L_H\MM$ with $(L_H\MM)_f\subset\mathrm{loc}_E(H)$.
\begin{proof} Suppose $\fromto{Y}{X}$ a fibration of $\mathrm{loc}_E(H)$, and suppose $\fromto{X'}{X}$ a morphism of $E$. To show that the pullback
\begin{equation*}
\xymatrix@C=18pt@R=18pt{Y'\ar[r]\ar[d]&Y\ar[d]\\
X'\ar[r]&X}
\end{equation*}
exists in $E$, it suffices by factorization (\ref{lem_Ertpropfunctfact}) to suppose that $\fromto{X'}{X}$ is a fibration as well. But then the pullback $Y'$ is a homotopy pullback, and since $\RMor(A,-)$ commutes with homotopy pullbacks, $Y'$ is also $H$-local.
\end{proof}
\end{prp}

\begin{cor}\label{cor_locEadmissiblehtypull} Admissible pullbacks of $\mathrm{loc}_E(H)$ are homotopy pullbacks.
\end{cor}

\begin{cor} If $\MM$ is right proper, then the quasifibrant objects form a right proper, admissibly left exact full subcategory $\mathrm{loc}(H)$ of $L_H\MM$.
\end{cor}

\begin{nul} Lastly, I now turn to a recognition principle for fibrations in $L_H\MM$ with codomains in $E$ or $\mathrm{loc}_{E}(H)$.
\end{nul}

\begin{prp}\label{prp:fibrationrecog} Suppose $p:\fromto{Y}{X}$ a fibration of $\MM$. For any fibrant replacement $p':\fromto{Y'}{X'}$ of $p$ (i.e., a morphism $p'$ between fibrant objects such that $p$ is isomorphic to $p'$ in $\Ho(\MM^1)$) consider the diagram
\begin{equation}\label{eqn:hopullback}
\xymatrix@C=18pt@R=18pt{
Y\ar[r]\ar[d]_p&Y'\ar[d]^{p'}\\
X\ar[r]&X'.
}
\end{equation}
\begin{enumerate}[(\ref{prp:fibrationrecog}.1)]
\addtocounter{enumi}{1}
\item If $X\in E$, then $p$ is a fibration of $L_H\MM$ if there exists a fibrant replacement $p':\fromto{Y'}{X'}$ of $p$ such that \eqref{eqn:hopullback} is a homotopy pullback square.
\item If $X\in\mathrm{loc}_E(H)$, then $p$ is a fibration of $L_H\MM$ if and only if for any fibrant replacement $p':\fromto{Y'}{X'}$, the diagram \eqref{eqn:hopullback} is a homotopy pullback square.
\end{enumerate}
\begin{proof} To prove the first assertion, factor $p'$ as a weak equivalence $\fromto{Y'}{Y''}$ of $L_H\MM$ followed by a fibration $\fromto{Y''}{X'}$ of $L_H\MM$; the weak equivalence $\fromto{Y'}{Y''}$ is even a weak equivalence of $\MM$ since it is a weak equivalence of local objects. Also factor $\fromto{X}{X'}$ as a weak equivalence $\fromto{X}{X''}$ in $E$ followed by a fibration $\fromto{X''}{X}$ in $E$. Pulling back $\fromto{Y''}{X'}$, we have the commutative diagram
\begin{equation*}
\xymatrix@C=18pt@R=18pt{
Y\ar[d]\ar[dr]\ar[rr]&&Y'\ar[d]\\
Z\ar[d]\ar[r]&Z''\ar[d]\ar[r]&Y''\ar[d]\\
X\ar[r]&X''\ar[r]&X'.
}
\end{equation*}
Since $X''$, $X'$, and $Y''$ are all objects of $E$, $Z''$ is an object of $E$ as well, and it follows from the right properness of $E$ that the weak equivalence $\fromto{X}{X''}$ is pulled back to a weak equivalence $\fromto{Z}{Z''}$ of $\MM$. The statement that \eqref{eqn:hopullback} is a homotopy pullback is equivalent to the assertion that the morphism $\fromto{Y}{Z''}$ is a weak equivalence of $\MM$. It now follows that the morphism $\fromto{Y}{Z}$ is a weak equivalence of $\MM$. Factor this map into a trivial cofibration $\fromto{Y}{Z'}$ of $L_H\MM$ followed by a fibration $\fromto{Z'}{Z}$ of $L_H\MM$; it follows that $\fromto{Y}{Z'}$ is in fact a trivial cofibration of $\MM$, and the retract argument thus completes proof.

The proof of the second statement begins similarly; factor $p'$ as before, and now factor the morphism $\fromto{X}{X'}$ as a weak equivalence $\fromto{X}{X''}$ in $\mathrm{loc}_E(H)$ followed by a fibration $\fromto{X''}{X}$ in $\mathrm{loc}_E(H)$ (which is in addition a weak equivalence in $L_H\MM$). Again pull back the fibration $\fromto{Y''}{X'}$:
\begin{equation*}
\xymatrix@C=18pt@R=18pt{
Y\ar[d]\ar[dr]\ar[rr]&&Y'\ar[d]\\
Z\ar[d]\ar[r]&Z''\ar[d]\ar[r]&Y''\ar[d]\\
X\ar[r]&X''\ar[r]&X'.
}
\end{equation*}
Now it follows from the right properness of $\mathrm{loc}_E(H)$ that the morphism $\fromto{Z''}{Y''}$ is a weak equivalence of $L_H\MM$, and it follows from the right properness of $E$ that $\fromto{Z}{Z''}$ is a weak equivalence of $\MM$. It thus follows that the morphism $\fromto{Y}{Z}$ is a weak equivalence of $L_H\MM$, and since $\fromto{Z}{X}$ and $p:\fromto{Y}{X}$ are each fibrations of $L_H\MM$, it follows that $\fromto{Y}{Z}$ is a weak equivalence of $\MM$; hence the composite morphism $\fromto{Y}{Z''}$ is a weak equivalence of $\MM$, so that \eqref{eqn:hopullback} is a homotopy pullback.
\end{proof}
\end{prp}

\subsection*{Application I: Presentations of combinatorial model categories} A presentation of a model $\XX$-category is a Quillen equivalence with a left Bousfield localization of a category of simplicial presheaves on an $\XX$-small category. A beautiful result of D. Dugger indicates that presentations exist for all tractable model categories. Hence any tracatable model category can (up to Quillen equivalence) be given a representation in terms of generators and relations. I recall Dugger's results here.

\begin{dfn} An \emph{$\XX$-presentation} $(K,H,F)$ of a model $\XX$-category $\MM$ consists of an $\XX$-small category $K$, an $\XX$-small set $H$ of homotopy classes of morphisms of $s\mathbf{Set}_{\XX}(K)$, and a left Quillen equivalence $F:\fromto{L_Hs\mathbf{Set}_{\XX}(K)_{\proj}}{\MM}$.
\end{dfn}

\begin{thm}[Dugger, \protect{\cite[Theorem 1.1]{MR1870516}}]\label{thm_duggerrep} Every $\XX$-combinatorial model $\XX$-category has an $\XX$-presentation.
\end{thm}

\begin{cor} An $\XX$-combinatorial model $\XX$-category has an $\XX$-small set of homotopy generators.
\begin{proof} By the theorem it is enough to show this for the projective model category of simplicial presheaves on an $\XX$-small category $C$; in this case, the images under the Yoneda embedding of the objects of $C$ provide such a set.
\end{proof}
\end{cor}

%\begin{nul} Using the existence of presentations, one may verify the following result. I include a brief sketch of a proof here; a more complete proof will appear elsewhere.
%\end{nul}

%\begin{cor} Suppose $\MM$ and $\NN$ model $\XX$-categories, and suppose $F:\fromto{L^H\MM}{L^H\NN}$ an equivalence of $s\mathbf{Set}_{\XX}$-categories. Then there exists a zigzag of Quillen equivalences
%\begin{equation*}
%\xymatrix@1@C=18pt{\MM\ar@<-0.5ex>[r]&\MM'\ar@<-0.5ex>[l]\ar@<0.5ex>[r]&\NN\ar@<0.5ex>[l]}
%\end{equation*}
%inducing $F$.
%\begin{proof}[Sketch of Proof] By the presentation theorem it suffices to give such a chain when $\MM=L_Hs\mathbf{Set}_{\XX}(K)_{\proj}$, $\NN=L_Js\mathbf{Set}_{\XX}(L)_{\proj}$, and $F$ is an equivalence $\fromto{\MM_{cf}}{\NN_{cf}}$ of $s\mathbf{Set}_{\XX}$-categories.

%Now consider the projective model structure on the categories $U(\MM_{cf})$ and $U(\NN_{cf})$ of small $s\mathbf{Set}_{\XX}$-functors $\fromto{\MM_{cf}}{s\mathbf{Set}_{\XX}}$ and $\fromto{\NN_{cf}}{s\mathbf{Set}_{\XX}}$, respectively, and form left Bousfield localizations $L_RU(\MM_{cf})$ and $L_SU(\NN_{cf})$ so that the fibrant objects are precisely the representables for each. Then the Yoneda embedding gives right Quillen equivalences $\fromto{\MM}{L_RU(\MM_{cf})}$ and $\fromto{\NN}{L_SU(\NN_{cf})}$. One now observes that $F$ induces an right Quillen equivalence $\fromto{L_RU(\MM_{cf})}{L_SU(\NN_{cf})}$.
%\end{proof}
%\end{cor}

\subsection*{Application II: Homotopy images} As a quirky demonstration of the usefulness of left Bousfield localizations, I offer the handy factorization result \ref{thm_factorcofhtymono}.

\begin{dfn}\label{dfn_htymonosubobj} Suppose $f:\fromto{X}{Y}$ a morphism of a model $\XX$-category $\MM$.
\begin{enumerate}[(\ref{dfn_htymonosubobj}.1)]\item The morphism $f$ is said to be a \emph{homotopy monomorphism} if the natural morphism $\fromto{X}{X\times^h_YX}$ is an isomorphism of $\Ho\MM$.
\item Dually, $f$ is a \emph{homotopy epimorphism} if the natural morphism $\fromto{Y\sqcup^{h,X}Y}{Y}$ is an isomorphism of $\Ho\MM$.
\item The \emph{homotopy image} of $f$ is a factorization of $f$ into a cofibration $\fromto{X}{f(X)}$ followed by a homotopy monomorphism $\fromto{f(X)}{Y}$ such that for any other such factorization $\xymatrix@1@C=10pt{X\ar[r]&X'\ar[r]&Y}$, there exists a unique morphism $\fromto{f(X)}{X'}$ in $\Ho\MM$.
\end{enumerate}
\end{dfn}

\begin{thm}\label{thm_factorcofhtymono} Suppose $\MM$ left proper and $\XX$-combinatorial; if $Y$ is a fibrant object of $\MM$, any morphism $f:\fromto{X}{Y}$ has a homotopy image.
\begin{proof} Let $G$ be an $\XX$-small set of homotopy generators for $\MM$. Write $G/Y$ for the disjoint union of the sets $\Mor_{\Ho\MM}(R,Y)$ over $R\in G$. Now write
\begin{equation*}
\nabla_{G/Y}:=\{\fromto{g\sqcup g}{g}\ |\ g\in G/Y\},
\end{equation*}
an $\XX$-small set. The $\MM$-model category $L_{\nabla_{G/Y}}(\MM/Y)$ then exists.

It now suffices to show that the fibrant objects of $L_{\nabla_{G/Y}}(\MM/Y)$ are precisely the fibrations $\fromto{X'}{Y}$ that are also homotopy monomorphisms, for if so, the homotopy image of a morphism $f:\fromto{X}{Y}$ is simply a fibrant replacement for $f$ in $L_{\nabla_{G/Y}}(\MM/Y)$.

Observe that the fibrant objects are precisely those fibrations $\fromto{X'}{Y}$ such that the natural morphism
\begin{equation*}
\fromto{\RMor_{(\MM/Y)}(Z,X')}{\star}
\end{equation*}
is a homotopy monomorphism of $\MM$ for any cofibrant object $Z$ of $\MM/Y$. Equivalently, a fibration $\fromto{X'}{Y}$ is fibrant in $L_{\nabla_{G/Y}}(\MM/Y)$ if and only if, for any object $Z$, the natural morphism
\begin{equation*}
\fromto{\RMor_{\MM}(Z,X')}{\RMor_{\MM}(Z,Y)}
\end{equation*}
is a homotopy monomorphism, whence the desired characterization of weak equivalences.
\end{proof}
\end{thm}

\subsection*{Application III: Resolutions of combinatorial model categories}{A surprisingly useful little trick is to replace a combinatorial model category $\MM$ with the category $\MM(K)$, equipped with a model structure that is Quillen equivalent to $\MM$ itself.}

\begin{nul} Suppose $\XX$ a universe, $\MM$ an $\XX$-combinatorial model category, and $K$ an $\XX$-small category.
\end{nul}

\begin{prp} There exists a model structure on $\MM(K)$ --- the \emph{$\colim$-resolution model structure} $\MM(K)_{\colim}$ --- such that the adjunction
\begin{equation*}
\colim:\xymatrix@1@C=18pt{\MM(K)_{\colim}\ar@<0.5ex>[r]&\MM\ar@<0.5ex>[l]}:\const
\end{equation*}
is a Quillen equivalence.
\begin{proof} Suppose $G$ an $\XX$-small set of cofibrant homotopy generators of $\MM(K)$, and let $R_{\MM}$ be a fibrant replacement functor for $\MM$. Set
\begin{equation*}
U_G:=\{\fromto{X}{\const R_{\MM}\colim X}\ |\ X\in G\},
\end{equation*}
and set $\MM(K)_{\colim}:=L_{U_G}\MM(K)_{\proj}$. Then $\colim:\fromto{\MM(K)}{\MM}$ factors through $\MM(K)_{\colim}$ since it sends elements of $U_G$ to weak equivalences. It is easy to see that the induced functor
\begin{equation*}
\RR\const:\fromto{\Ho\MM}{\Ho\MM(K)_{\colim}}
\end{equation*}
is fully faithful and essentially surjective.
\end{proof}
\end{prp}

\subsection*{Application IV: Homotopy limits of right Quillen presheaves}{The category of right sections of a right Quillen presheaf $\FF$ with its projective model structure is to be thought of as the $(\infty,1)$-categorical lax limit of $\FF$. The $(\infty,1)$-categorical limit --- or homotopy limit --- of $\FF$ is a left Bousfield localization of the category of right sections.}

\begin{nul} Suppose $\XX$ a universe. Suppose $K$ and $\XX$-small category, and suppose $\FF$ an $\XX$-combinatorial right Quillen presheaf on $K$.
\end{nul}

\begin{dfn} A right section $(X,\phi)$ of $\FF$ is said to be \emph{homotopy cartesian} if for any morphism $f:\fromto{\ell}{k}$ of $K$, the morphism
\begin{equation*}
\phi_f^h:\fromto{X_{\ell}}{\RR f^{\star}X_k}
\end{equation*}
is an isomorphism of $\Ho\FF_{\ell}$.
\end{dfn}

\begin{thm}\label{prp:holimmodcat} There exists an $\XX$-combinatorial model structure on the category $\Sect^R\FF$ --- the homotopy limit structure $\Sect_{\holim}^R\FF$ --- satisfying the following conditions.
\begin{enumerate}[(\ref{prp:holimmodcat}.1)]
\item\addtocounter{equation}{1} The cofibrations are exactly the projective cofibrations.
\item\addtocounter{equation}{1} The fibrant objects are the projective fibrant right sections that are homotopy cartesian.
\item\addtocounter{equation}{1} The weak equivalences between fibrant objects are precisely the objectwise weak equivalences.
\end{enumerate}
\begin{proof} For every $k\in\Obj K$, let $G_k$ be an $\XX$-small set of cofibrant homotopy generators of $\FF_k$. For each object $k$ of $K$, there is a Quillen adjunction
\begin{equation*}
D_k:\xymatrix@1@C=18pt{\FF_k\ar@<0.5ex>[r]&\Sect^R\FF\ar@<0.5ex>[l]}:E_k
\end{equation*}
where $E_k(X,\phi)=X_k$; one verifies that for any object $A$ of $\FF_{\ell}$,
\begin{equation*}
E_kD_{\ell}A\cong\Coprod_{f:\ell\to k}f_!A.
\end{equation*}
Hence for every $f:\fromto{\ell}{k}$ and any object $A$ of $\FF_{\ell}$, there is a canonical morphism $\fromto{f_!A}{E_kD_{\ell}A}$, and, by adjunction, a canonical morphism $r_{f,A}:\fromto{D_k(f_!A)}{D_{\ell}A}$.

Now define the $\XX$-small set
\begin{equation*}
H:=\{r_{f,A}:\fromto{D_k(f_!A)}{D_{\ell}A}\ |\ [f:\fromto{\ell}{k}]\in K, A\in G_{\ell}\}.
\end{equation*}
The claim is that $L_H\Sect^R_{\proj}\FF$ is the model category of the theorem.

To verify this claim, it suffices to check that the fibrant objects are as described. Indeed, a right section $X$ is fibrant if and only if it is fibrant in $\Sect^R_{\proj}\FF$, and, for any morphism $f:\fromto{\ell}{k}$ of $K$ and any $A\in G_{\ell}$,
\begin{equation*}
\fromto{\RMor_{\Sect^R_{\proj}\FF}(\LL D_{\ell}A,X)}{\RMor_{\Sect^R_{\proj}\FF}(\LL D_k\LL f_!A,X)},
\end{equation*}
or equivalently
\begin{equation}\label{eqn:genrtglue}
\fromto{\RMor_{\FF_k}(A,\RR E_{\ell}X)}{\RMor_{\FF_k}(A,\RR f^{\star}\RR E_kX)}
\end{equation}
is an isomorphism of $\Ho s\mathbf{Set}_{\XX}$. Since the elements of $G_{\ell}$ generate $\FF_{\ell}$ by homotopy colimits, it follows immediately that \eqref{eqn:genrtglue} is an isomorphism of $\Ho s\mathbf{Set}_{\XX}$ for any object $A$ of $F_{\ell}$, whence it follows that $X$ is homotopy cartesian, as desired.
\end{proof}
\end{thm}

\begin{exm} Denote by $\NN$ the category whose objects are nonnegative integers, in which there is a unique morphism $\fromto{m}{n}$ if and only if $m\leq n$. Consider the right Quillen presheaf
\begin{equation*}
\mathbf{\Omega}:\xymatrix@1@C=18pt@R=0pt{\NN^{\op}\ar[r]&\mathbf{Cat}_{\YY}\\
n\,\ar@{|->}[r]&(\star/s\mathbf{Set}_{\XX})\\
[n\leq m]\,\ar@{|->}[r]&\Omega^{m-n},
}
\end{equation*}
where of course $\Omega^{m-n}:=\MOR_{(\star/s\mathbf{Set})}((S^1)^{\wedge(m-n)},-)$. One verifies easily that $\Sect_{\holim}^R\mathbf{\Omega}$ is simply the usual model category of spectra.
\end{exm}

\begin{nul} Note that the results of this section say nothing about the homotopy limits of left Quillen presheaves. As an $(\infty,1)$-category, such a homotopy limit should be a corefexive sub-$(\infty,1)$-category of the $(\infty,1)$-categorical lax limit; hence it is more properly modeled as a right Bousfield localization. This is a somewhat delicate issue, which we will take up elsewhere.
\end{nul}

\subsection*{Application V: Postnikov towers for simplicial model categories}{Any fibrant simplicial set $X$ has a Postnikov tower
\begin{equation*}
\xymatrix@C=18pt{\dots\ar[r]&X\langle n\rangle\ar[r]&X\langle n-1\rangle\ar[r]&\dots\ar[r]&X\langle 1\rangle\ar[r]&X\langle 0\rangle},
\end{equation*}
in which each $X\langle n\rangle$ is an $n$-type, and for any $0\leq j<n$, the morphism $f:\fromto{X\langle n\rangle}{X\langle n-1\rangle}$ induces an isomoprhism
\begin{equation*}
\fromto{\pi_j(X\langle n\rangle,x)}{\pi_j(X\langle n-1\rangle,f(x))}
\end{equation*}
for any point $x$ of $X$. An analogous construction can be made in any combinatorial, left proper, simplicial model category.}

\begin{nul} Suppose $\XX$ a universe, $\MM$ an $\XX$-combinatorial, left proper, simplicial model $\XX$-category.
\end{nul}

\begin{dfn} For any integer $n\geq -1$, an object $X$ of $\MM$ is \emph{$n$-truncated} if for any object $Z$ of $\MM$, the simplicial set $\RMor_{\MM}(Z,X)$ is an $n$-type; the object $X$ is \emph{truncated} if and only if it is $n$-truncated for some $n$.
\end{dfn}

\begin{prp}\label{prp:ntruncated} For any integer $n\geq -1$, there exists a combinatorial, left proper, simplicial model structure on the category $\MM$ --- the \emph{$n$-truncated model structure} $\MM_{\leq n}$ --- satisfying the following conditions.
\begin{enumerate}[(\ref{prp:ntruncated}.1)]
\item The cofibrations of $\MM_{\leq n}$ are precisely the cofibrations of $\MM$.
\item The fibrant objects of $\MM_{\leq n}$ are precisely the fibrant, $n$-truncated objects of $\MM$.
\item The weak equivalences between the fibrant objects are precisely the weak equivalences of $\MM$.
\end{enumerate}
\begin{proof} Let $G$ be an $\XX$-small set of cofibrant homotopy generators of $\MM$. Then one verifies easily that the $n$-truncated model structure is the left Bousfield localization $L_{H(n)}\MM$ with respect to the set
\begin{equation*}
H(n):=\{\fromto{S^j\otimes X}{X}\ |\ X\in G, 0\leq n<j\}.\qedhere
\end{equation*}
\end{proof}
\end{prp}

\begin{prp}\label{prp:postnikov} There exists a combinatorial, left proper, simplicial model structure on the presheaf category $\MM(\NN)$ --- the \emph{Postnikov model structure} $\MM(\NN)_{\Post}$ --- satisfying the following conditions.
\begin{enumerate}[(\ref{prp:postnikov}.1)]
\item The cofibrations of $\MM(\NN)_{\Post}$ are precisely the objectwise cofibrations.
\item\label{item:postnikovfibrant} The fibrant objects are those sequences
\begin{equation*}
\xymatrix@C=18pt{\dots\ar[r]&X(n)\ar[r]&X(n-1)\ar[r]&\dots\ar[r]&X(1)\ar[r]&X(0)}
\end{equation*}
satisfying the following conditions.
\begin{enumerate}[(\ref{prp:postnikov}.\ref{item:postnikovfibrant}.1)]
\item The object $X(0)$ is fibrant.
\item Each morphism $\fromto{X(n)}{X(n-1)}$ is a fibration.
\item Each object $X(n)$ is \emph{$n$-truncated}, in the sense that for any object $Z$ of $\MM$, $\RMor_{\MM}(Z,X(n))$ is an $n$-type.
\item For any object $Z$ of $\MM$ and any integers $0\leq j<n$, the morphism $\sigma:\fromto{X\langle n\rangle}{X\langle n-1\rangle}$ induces an isomorphism
\begin{equation*}
\fromto{\pi_j(\RMor_{\MM}(Z,X(n)),\phi)}{\pi_j(\RMor_{\MM}(Z,X(n-1)),\sigma\circ\phi)}.
\end{equation*}
\end{enumerate}
\item The weak equivalences between the fibrant objects are precisely the objectwise weak equivalences.
\end{enumerate}
\begin{proof} Let $G$ be an $\XX$-small set of cofibrant homotopy generators of $\MM$. Then I claim that the Postnikov model structure is the left Bousfield localization $L_H\MM(\NN)_{\inj}$ of the injective model structure with respect to the set
\begin{equation*}
H:=\{\fromto{R^n(S^j\otimes X)}{R^nX}, \fromto{R^{j-1}(S^n\otimes X)}{R^j(S^n\otimes X)}\ |\ X\in G, 0\leq n<j\},
\end{equation*}
where $R^n:\fromto{\MM}{\MM(\NN)}$ is the left adjoint to the evaluation-at-$n$ functor $\goesto{X}{X(n)}$:
\begin{equation*}
R^nX(m):=\begin{cases} \varnothing&\textrm{if }m>n\\
X&\textrm{if }m\leq n.
\end{cases}
\end{equation*}

It suffices to verify that the fibrant objects are as described, and this follows directly from an elementary adjunction exercise.
\end{proof}
\end{prp}

\begin{cor} For any integer $n\geq -1$, the evaluation-at-$n$ functor
\begin{equation*}
\xymatrix@1@R=0pt@C=18pt{\MM(\NN)_{\Post}\ar[r]&\MM_{\leq n}\\
A\,\ar@{|->}[r]&A(n)}
\end{equation*}
is left Quillen.
\end{cor}

\begin{dfn}\label{dfn:Posttowerhypercomplete} Suppose $\MM$ a combinatorial, left proper, simplicial model category.
\begin{enumerate}[(\ref{dfn:Posttowerhypercomplete}.1)]
\item The \emph{Postnikov tower} of an object $X$ in $\MM$ is a fibrant model for
\begin{equation*}
\const X:=[\xymatrix@C=18pt{\dots\ar@{=}[r]&X\ar@{=}[r]&X\ar@{=}[r]&\dots\ar@{=}[r]&X\ar@{=}[r]&X}]
\end{equation*}
in $\MM(\NN)_{\Post}$; it will be denoted
\begin{equation*}
\xymatrix@C=18pt{\dots\ar[r]&X\langle n\rangle\ar[r]&X\langle n-1\rangle\ar[r]&\dots\ar[r]&X\langle 1\rangle\ar[r]&X\langle 0\rangle}
\end{equation*}
\item One says that a morphism $\fromto{X}{Y}$ of $\MM$ is \emph{$\infty$-connected} if it induces an equivalence $\fromto{\const X}{\const Y}$ in $\MM(\NN)_{\Post}$.
\item One says that $\MM$ is \emph{hypercomplete} if every object is the homotopy limit of its Postnikov tower, i.e., if $\holim\circ\LL\const$ is isomorphic to the identity functor on $\Ho\MM$.
\end{enumerate}
\end{dfn}

\begin{prp} There exists a \emph{hypercompletion} $\fromto{\MM}{\MM^{\wedge}}$ of $\MM$, which is the initial object in the category of left Quillen functors $\fromto{\MM}{\PP}$ under which $\infty$-connected morphisms are sent to weak equivalences.
\begin{proof} Simply define $\MM^{\wedge}$ as the left Bousfield localization of $\MM$ with respect to the set $\{\fromto{\lim R\const X}{X}\ |\ X\in G\}$, where $G$ is again an $\XX$-small set of cofibrant homotopy generators of $\MM$, and $R$ is a fibrant replacement functor.
\end{proof}
\end{prp}

\begin{lem}\label{lem:hypercmpltequiv} The following conditions are equivalent.
\begin{enumerate}[(\ref{lem:hypercmpltequiv}.1)]
\item\label{item:MMhypercmplt} The model category $\MM$ is hypercomplete.
\item\label{item:MMhypercmpltion} The hypercompletion functor $\fromto{\MM}{\MM^{\wedge}}$ is a Quillen equivalence.
\item\label{item:inftyconnwe} Every $\infty$-connected morphism of $\MM$ is a weak equivalence.
\item\label{item:wedetectbytruncs} A morphism $\fromto{Z}{W}$ of $\MM$ is a weak equivalence if the induced morphism $\fromto{\Mor_{\Ho\MM}(W,X)}{\Mor_{\Ho\MM}(Z,X)}$ is a bijection for every truncated object $X$ of $\MM$.
\end{enumerate}
\begin{proof} It is readily apparent that (\ref{lem:hypercmpltequiv}.\ref{item:MMhypercmplt}) and (\ref{lem:hypercmpltequiv}.\ref{item:MMhypercmpltion}) are equivalent conditions.

If $\MM$ is hypercomplete, then every $\infty$-connected morphism of $\MM$ is a homotopy limit of weak equivalences; hence (\ref{lem:hypercmpltequiv}.\ref{item:MMhypercmplt}) implies (\ref{lem:hypercmpltequiv}.\ref{item:inftyconnwe}). Conversely, one sees that the morphism $\fromto{X}{\lim_{n\in\NN}X\langle n\rangle}$ is automatically $\infty$-connected, so in fact (\ref{lem:hypercmpltequiv}.\ref{item:MMhypercmplt}) and (\ref{lem:hypercmpltequiv}.\ref{item:inftyconnwe}) are equivalent conditions.

For any morphism $\fromto{Z}{W}$, if $\fromto{\Mor_{\Ho\MM}(W,X)}{\Mor_{\Ho\MM}(Z,X)}$ is a bijection for every truncated object $X$, then the induced morphism $\fromto{\const Z}{\const W}$ in $\MM(\NN)_{\Post}$ is easily seen to be a weak equivalence, whence it follows that (\ref{lem:hypercmpltequiv}.\ref{item:inftyconnwe}) implies (\ref{lem:hypercmpltequiv}.\ref{item:wedetectbytruncs}).

Suppose $Z$ an object of $\MM$. To show that (\ref{lem:hypercmpltequiv}.\ref{item:wedetectbytruncs}) implies (\ref{lem:hypercmpltequiv}.\ref{item:inftyconnwe}), it suffices to show that for any truncated object $X$ of $\MM$ and any object $Z$ of $\MM$,
\begin{equation*}
\fromto{\Mor_{\Ho\MM}(\holim_{n\in\NN^{\op}}Z\langle n\rangle,X)}{\Mor_{\Ho\MM}(Z,X)}
\end{equation*}
is a bijection. For this, one sees that if $X$ is $k$-truncated, then one may factor the above morphism as a series of bijections
\begin{eqnarray}
\Mor_{\Ho\MM}(\holim_{n\in\NN^{\op}}Z\langle n\rangle,X)&\cong&\Mor_{\Ho\MM}((\holim_{n\in\NN^{\op}}Z\langle n\rangle)\langle k\rangle,X)\nonumber\\
&\cong&\Mor_{\Ho\MM}(Z\langle k\rangle,X)\nonumber\\
&\cong&\Mor_{\Ho\MM}(Z,X).\nonumber
\end{eqnarray}

This completes the proof.
\end{proof}
\end{lem}

\section{Enriched model categories and enriched left Bousfield localization}

\subsection*{Symmetric monoidal model categories and enrichments} I present a thorough, if somewhat terse, review of the basic theory of symmetric monoidal and enriched model categories.

\begin{nul} Again suppose $\XX$ a universe.
\end{nul}

\begin{dfn}\label{dfn_symmonmod}
\begin{enumerate}[(\ref{dfn_symmonmod}.1)]
\item\label{Quillenadjof2vars} Suppose $\DD$, $\EE$, and $\FF$ model $\XX$-categories. Suppose
\begin{equation*}
\otimes:\fromto{\DD\times\EE}{\FF}\qquad\MOR:\fromto{\EE^{\op}\times\FF}{\DD}\qquad\mor:\fromto{\DD^{\op}\times\FF}{\EE}
\end{equation*}
form an adjunction of two variables. Then $(\otimes,\MOR,\mor)$ is a \emph{Quillen adjunction of two variables} if the following axiom holds.
\begin{enumerate}[(\ref{dfn_symmonmod}.\ref{Quillenadjof2vars}.1)]
\item\label{pushoutprod} (Pushout-product axiom) For any pair of cofibrations $c:\fromto{Q}{R}$ of $\DD$ and $d:\fromto{S}{T}$ of $\EE$, the \emph{pushout-product}
\begin{equation*}
c\Box d:\fromto{(Q\otimes T)\sqcup^{(Q\otimes S)}(R\otimes S)}{R\otimes T}
\end{equation*}
is a cofibration of $\FF$ that is trivial if either $f$ or $g$ is.\footnote{If $S$ and $T$ are sets of morphisms, it will be convenient to denote by $S\Box T$ the set of morphisms of the form $f\Box g$ for $f\in S$ and $g\in T$.}
\end{enumerate}
\item\label{symmonmodcat} A \emph{symmetric monoidal model $\XX$-category} is a symmetric monoidal closed $\XX$-category $(\VV,\otimes_{\VV},\mathbf{1}_{\VV},\MOR_{\VV})$ \cite[Definitions 6.1.1--3]{MR1313497}, equipped with a model structure such that the following axioms hold.
\begin{enumerate}[(\ref{dfn_symmonmod}.\ref{symmonmodcat}.1)]
\item The tuple $(\otimes_{\VV},\MOR_{\VV},\MOR_{\VV})$ is a Quillen adjunction of two variables.
\item (Unit axiom) For any object $A$, the canonical morphism
\begin{equation*}
\xymatrix@1@C=18pt{Q_{\VV}\mathbf{1}_{\VV}\otimes_{\VV}A\ar[r]&\mathbf{1}_{\VV}\otimes_{\VV}A\ar[r]&A}
\end{equation*}
is a weak equivalence for some cofibrant replacement $\fromto{Q_{\VV}\mathbf{1}_{\VV}}{\mathbf{1}_{\VV}}$.
\end{enumerate}
\item An \emph{internal closed model category} is a symmetric monoidal model category $\VV$ in which $\otimes_{\VV}$ is the categorical (cartesian) product $\times$.
\item\label{modVcat} If $\VV$ is a symmetric monoidal model $\XX$-category, then a \emph{model $\VV$-category} is a tensored and cotensored $\VV$-category $(\CC,\MOR_{\CC}^{\VV},\otimes_{\CC}^{\VV},\mor_{\CC}^{\VV})$ \cite[Definitions 6.2.1, 6.5.1]{MR1313497}, equipped with a model structure on the underlying $\XX$-category of $\CC$ such that the following axioms hold.
\begin{enumerate}[{(\ref{dfn_symmonmod}.\ref{modVcat}.}1{)}]
\item The tuple $(\otimes_{\CC}^{\VV},\MOR_{\CC}^{\VV},\mor_{\CC}^{\VV})$ is a Quillen adjunction of two variables.
\item (Unit axiom) For any object $X$ of $\CC$, the canonical morphism
\begin{equation*}
\xymatrix@1@C=18pt{Q_{\VV}\mathbf{1}_{\VV}\otimes_{\CC}^{\VV}X\ar[r]&\mathbf{1}_{\VV}\otimes_{\CC}^{\VV}X\ar[r]&X}
\end{equation*}
is a weak equivalence for some cofibrant replacement $\fromto{Q_{\VV}\mathbf{1}_{\VV}}{\mathbf{1}_{\VV}}$.
\end{enumerate}
\item A \emph{simplicial model $\XX$-category} $(\MM,\Map_{\MM},\otimes_{\MM},\mor_{\MM})$ is a model $s\mathrm{Set}_{\XX}$-category.
\item\label{QuillenVfunctor}Suppose $\VV$ a symmetric monoidal model category and $\CC$ and $\CC'$ model $\VV$-categories.
\begin{enumerate}[(\ref{dfn_symmonmod}.\ref{QuillenVfunctor}.1)]
\item A left (respectively, right) $\VV$-adjoint $F:\fromto{\CC}{\CC'}$ \cite[Definition 6.7.1]{MR1313497} whose underlying functor $F_0$ is a left (resp., right) Quillen functor will be called a \emph{left} (resp., \emph{right}) \emph{Quillen $\VV$-functor}.
\item If $F_0$ is in addition a Quillen equivalence, then $F$ will be called a \emph{left} (resp., \emph{right}) \emph{Quillen $\VV$-equivalence}.
\end{enumerate}
\end{enumerate}
\end{dfn}

\begin{ntn} Of course the sub- and superscripts on $\otimes$, $\mathbf{1}$, and $\MOR$ will be dropped if no confusion can result, and by the standard harmless abuse, we will refer to $\VV$ alone as the symmetric monoidal model category.
\end{ntn}

\begin{lem}[\protect{\cite[Lemma 4.2.2]{MR99h:55031}}]\label{lem_tfaeQ2vars} The following are equivalent for three model $\XX$-categories $\DD$, $\EE$, and $\FF$ and an adjunction of two variables
\begin{equation*}
\otimes:\fromto{\DD\times\EE}{\FF}\qquad\MOR:\fromto{\EE^{\op}\times\FF}{\DD}\qquad\mor:\fromto{\DD^{\op}\times\FF}{\EE}.
\end{equation*}
\begin{enumerate}[(\ref{lem_tfaeQ2vars}.1)]
\item The tuple $(\otimes,\MOR,\mor)$ is a Quillen adjunction of two variables.
\item For any cofibration $d:\fromto{S}{T}$ of $\EE$ and any fibration $b:\fromto{V}{U}$ of $\FF$, the morphism
\begin{equation*}
\MOR_{\Box}(d,b):\fromto{\MOR(T,V)}{\MOR(S,V)\times_{\MOR(S,U)}\MOR(T,U)}
\end{equation*}
is a fibration that is trivial if either $d$ or $b$ is.
\item For any cofibration $c:\fromto{Q}{R}$ of $\DD$ and any fibration $b:\fromto{V}{U}$ of $\FF$, the morphism
\begin{equation*}
\mor_{\Box}(c,b):\fromto{\mor(R,V)}{\mor(Q,V)\times_{\mor(Q,U)}\mor(R,U)}
\end{equation*}
is a fibration that is trivial if either $c$ or $b$ is.
\end{enumerate}
\end{lem}

\begin{lem}[\protect{\cite[Lemma 4.2.4]{MR99h:55031}}]\label{lem_Q2varsgencof} Suppose $\DD$ an $\XX$-cofibrantly generated model $\XX$-category, with generating cofibrations $I_{\DD}$ and generating trivial cofibrations $J_{\DD}$. Suppose $\EE$ an $\XX$-cofibrantly generated model $\XX$-category, with generating cofibrations $I_{\EE}$ and generating trivial cofibrations $J_{\EE}$. Suppose $\FF$ a model $\XX$-category, and suppose
\begin{equation*}
\otimes:\fromto{\DD\times\EE}{\FF}\qquad\MOR:\fromto{\EE^{\op}\times\FF}{\DD}\qquad\mor:\fromto{\DD^{\op}\times\FF}{\EE}
\end{equation*}
form an adjunction of two variables. Then $(\otimes,\MOR,\mor)$ is a \emph{Quillen adjunction of two variables} if and only if the following conditions hold.
\begin{enumerate}[(\ref{lem_Q2varsgencof}.A)]
\item $I_{\DD}\Box I_{\EE}$ consists only of cofibrations.
\item $I_{\DD}\Box J_{\EE}$ consists only of weak equivalences.
\item $J_{\DD}\Box I_{\EE}$ consists only of weak equivalences.
\end{enumerate}
\end{lem}

\begin{lem}\label{lem_symmonleftQuillen} Suppose $\VV$ and $\CC$ symmetric monoidal model $\XX$-categories, wherein the unit $1_{\CC}$ is cofibrant. Then a model $\VV$-category structure on $\CC$ is equivalent to a Quillen adjunction
\begin{equation*}
\mathrm{real}:\xymatrix@1@C=18pt{\VV\ar@<0.5ex>[r]&\CC\ar@<0.5ex>[l]}:\Pi
\end{equation*}
in which the left adjoint $\mathrm{real}$ is symmetric monoidal.
\begin{proof} Suppose $(\otimes_{\CC}^{\VV},\MOR_{\CC}^{\VV},\mor_{\CC}^{\VV})$ a model $\VV$-category structure on $\CC$. Set
\begin{equation*}
\mathrm{real}:=-\otimes_{\CC}^{\VV}\mathbf{1}_{\CC}\qquad\Pi:=\MOR_{\CC}^{\VV}(\mathbf{1}_{\CC},-).
\end{equation*}
For any objects $K$ and $L$ of $\VV$, one verifies easily that the objects $(K\otimes_{\VV}L)\otimes_{\CC}^{\VV}\mathbf{1}_{\CC}$ and $(K\otimes_{\CC}^{\VV}\mathbf{1}_{\CC})\otimes_{\VV}(L\otimes_{\CC}^{\VV}\mathbf{1}_{\CC})$ corepresent the same functor. Since $\mathbf{1}_{\CC}$ is cofibrant, the pushout-product axiom implies that $\mathrm{real}$ is left Quillen.

On the other hand, suppose
\begin{equation*}
\mathrm{real}:\xymatrix@1@C=18pt{\VV\ar@<0.5ex>[r]&\CC\ar@<0.5ex>[l]}:\Pi
\end{equation*}
a Quillen adjunction in which the left adjoint $\mathrm{real}$ is symmetric monoidal. For $K$ and object of $\VV$, and $X$ and $Y$ objects of $\CC$, set
\begin{eqnarray}
\Mor_{\CC}^{\VV}(X,Y)&:=&\Pi\MOR_{\CC}(X,Y),\nonumber\\
K\otimes_{\CC}^{\VV}X&:=&\mathrm{real}(K)\otimes_{\CC}X,\nonumber\\
\mor_{\CC}^{\VV}(K,Y)&:=&\MOR_{\CC}(\mathrm{real}(K),Y).\nonumber
\end{eqnarray}
The pushout-product axiom for $\otimes_{\CC}$ implies the pushout-product axiom for $\otimes_{\CC}^{\VV}$.

These two definitions are inverse to one another.
\end{proof}
\end{lem}

\begin{lem} Suppose $\CC$ a simplicial, internal model $\XX$-category. Then for any $\XX$-small set $S$, the object $\mathrm{real}(S)$ is canonically isomorphic to the copower $S\cdot\star$.
\begin{proof} The functor $\mathrm{real}$ is a left adjoint and thus respects copowers; it is symmetric monoidal and thus preserves terminal objects.
\end{proof}
\end{lem}

\begin{prp} Suppose $\VV$ internal. Then for any object $Y$ of $\VV$, the comma category $\VV/Y$ is a $\VV$-model category with
\begin{eqnarray}
Z\otimes_{(\VV/Y)}^{\VV}X&:=&Z\times X\nonumber\\
\mor_{(\VV/Y)}^{\VV}(Z,X')&:=&Y\times_{\MOR_{\VV}(Z,Y)}\MOR_{\VV}(Z,X')\nonumber\\
\MOR_{(\VV/Y)}^{\VV}(X,X')&:=&\star\times_{\MOR_{\VV}(X,Y)}\MOR_{\VV}(X,X')\nonumber
\end{eqnarray}
for any object $Z$ of $\VV$ and any objects $X$ and $X'$ of $\VV/Y$.
\begin{proof} A morphism of the comma category $\VV/Y$ is a cofibration, fibration, or weak equivalence if and only if its image under the forgetful functor $\fromto{(\VV/Y)}{\VV}$ is so; hence the pushout-product and unit axioms for $(\VV/Y)$ follow directly from those for $\VV$.
\end{proof}
\end{prp}

\begin{nul} I will be interested in localizing enriched model categories using derived mapping objects in lieu of the homotopy function complexes. The result has a universal property that is rather different from the ordinary left Bousfield localization. Suppose $(\VV,\otimes_{\VV},\mathbf{1}_{\VV},\MOR_{\VV})$ a symmetric monoidal model $\XX$-category, and suppose $\CC$ a model $\VV$-category.
\end{nul}

\begin{dfn}\label{dfn_RHOM}
\begin{enumerate}[(\ref{dfn_RHOM}.1)]
\item The left Kan extension (if it exists) of the composite
\begin{equation*}
\xymatrix@1@C=18pt{\CC^{\op}\times\CC\ar[r]&\VV\ar[r]&\Ho\VV}
\end{equation*}
along the localization functor $\fromto{\CC^{\op}\times\CC}{\Ho(\CC^{\op}\times\CC)}$ is the \emph{derived mapping object functor}, denoted $\RMOR_{\CC}^{\VV}$.
\item Suppose $Q$ and $R$ cofibrant and fibrant replacement functors for $\CC$; then the $(Q,R)$-\emph{derived mapping object functor} is the functor
\begin{equation*}
\xymatrix@1@C=18pt@R=0pt{\RMOR_{\CC,Q,R}:\CC^{\op}\times\CC\ar[r]&\Ho\VV\\
\qquad\qquad\qquad(A,B)\,\ar@{|->}[r]&\MOR_{\CC}(QA,RB).}
\end{equation*}
\end{enumerate}
\end{dfn}

\begin{prp} Suppose $Q$ and $R$ cofibrant and fibrant replacement functors for $\CC$, respectively. Then the functor $\RMOR_{\CC}$ exists, and, up to isomorphism of functors, $\RMOR_{\CC,Q,R}$ factors through it.
\begin{proof} This is an immediate consequence of \cite[Proposition 8.4.8]{MR2003j:18018}.
\end{proof}
\end{prp}

\begin{lem}\label{lem_RMORdetectsws} The following are equivalent for a morphism $\fromto{A}{B}$ of $\CC$.
\begin{enumerate}[(\ref{lem_RMORdetectsws}.1)]
\item\label{item_ABwe} The morphism $\fromto{A}{B}$ is a weak equivalence.
\item\label{item_RMORweright} For any fibrant object $Z$ of $\CC$, the induced morphism
\begin{equation*}
\fromto{\RMOR_{\CC}(B,Z)}{\RMOR_{\CC}(A,Z)}
\end{equation*}
is an isomorphism of $\Ho\VV$.
\item\label{item_RMORweleft} For any cofibrant object $X$ of $\CC$, the induced morphism
\begin{equation*}
\fromto{\RMOR_{\CC}(X,A)}{\RMOR_{\CC}(X,B)}
\end{equation*}
is an isomorphism of $\Ho\VV$.
\end{enumerate}
\begin{proof} That (\ref{lem_RMORdetectsws}.\ref{item_ABwe}) implies both (\ref{lem_RMORdetectsws}.\ref{item_RMORweright}) and (\ref{lem_RMORdetectsws}.\ref{item_RMORweleft}) follows from the pushout-product axiom.

To show that (\ref{lem_RMORdetectsws}.\ref{item_ABwe}) follows from (\ref{lem_RMORdetectsws}.\ref{item_RMORweright}), suppose $\fromto{A}{B}$ a morphism of $\CC$ such that for any fibrant object $Z$ of $\CC$, the induced morphism
\begin{equation*}
\fromto{\RMOR_{\CC}(B,Z)}{\RMOR_{\CC}(A,Z)}
\end{equation*}
is an isomorphism of $\Ho\VV$; one may clearly assume that $A$ and $B$ are cofibrant, so that the morphism
\begin{equation*}
\fromto{\MOR_{\CC}(B,Z)}{\MOR_{\CC}(A,Z)}
\end{equation*}
is a weak equivalence of $\VV$. Applying $\RMor_{\VV}(1_{\VV},-)$ yields an isomorphism
\begin{equation*}
\fromto{\RMor_{\CC}(B,Z)}{\RMor_{\CC}(A,Z)}
\end{equation*}
of $\Ho s\mathrm{Set}_{\XX}$ for any fibrant object $Z$ of $\CC$. Now apply \cite[17.7.7]{MR2003j:18018}.

The very same argument, mutatis mutandis, shows that (\ref{lem_RMORdetectsws}.\ref{item_ABwe}) follows from (\ref{lem_RMORdetectsws}.\ref{item_RMORweleft}).
\end{proof}
\end{lem}

\subsection*{The enriched left Bousfield localization}{Here I define enriched Bousfield localizations, and I prove an existence theorem.}

\begin{nul} Suppose $\XX$ a universe, $\VV$ a symmetric monoidal model $\XX$-category, and $\CC$ a model $\VV$-category.
\end{nul}

\begin{dfn} Suppose $H$ a set of homotopy classes of morphisms of $\CC$. A \emph{left Bousfield localization of $\CC$ with respect to $H$ enriched over $\VV$} is a model $\VV$-category $L_{(H/\VV)}\CC$, equipped with a left Quillen $\VV$-functor $\fromto{\CC}{L_{(H/\VV)}\CC}$ that is initial among left Quillen $\VV$-functors $F:\fromto{\CC}{\DD}$ to model $\VV$-categories $\DD$ such that for any $f$ representing a class in $H$, $Ff$ is a weak equivalence in $\DD$.
\end{dfn}

\begin{lem} When it exists, an enriched left Bousfield localization is unique up to a unique $\VV$-isomorphism of model $\VV$-categories under $\CC$.
\begin{proof} Initial objects are essentially unique.
\end{proof}
\end{lem}

\begin{nul} Suppose, for the remainder of this section, $H$ a set of homotopy classes of morphisms of $\CC$.
\end{nul}

\begin{dfn}\label{dfn_HVlocal}
\begin{enumerate}[(\ref{dfn_HVlocal}.1)]
\item An object $Z$ of $\CC$ is \emph{$(H/\VV)$-local} if for any morphism $\fromto{A}{B}$ representing an element of $H$, the morphism
\begin{equation*}
\fromto{\RMOR_{\CC}(B,Z)}{\RMOR_{\CC}(A,Z)}
\end{equation*}
is an isomorphism of $\Ho\VV$.
\item A morphism $\fromto{A}{B}$ of $\CC$ is an \emph{$(H/\VV)$-local equivalence} if for any fibrant $(H/\VV)$-local object $Z$, the morphism
\begin{equation*}
\fromto{\RMOR_{\CC}(B,Z)}{\RMOR_{\CC}(A,Z)}
\end{equation*}
is an isomorphism of $\Ho\VV$.
\end{enumerate}
\end{dfn}

\begin{thm}\label{thm_elbexist} Suppose that following conditions are satisfied.
\begin{enumerate}[(\ref{thm_elbexist}.A)]
\item\label{CClpropXcomb} The model $\VV$-category $\CC$ is left proper and $\XX$-tractable.
\item\label{HXXsmall} The set $H$ is $\XX$-small.
\item\label{VVcofgen} The model $\XX$-category $\VV$ is $\XX$-tractable.
\end{enumerate}
Then the enriched left Bousfield localization $L_{(H/\VV)}\CC$ exists, and satisfies the following conditions.
\begin{enumerate}[(\ref{thm_elbexist}.1)]
\item The model category $L_{(H/\VV)}\CC$ is left proper, $\XX$-combinatorial, and right tractable.
\item As a $\VV$-category, $L_{(H/\VV)}\CC$ is simply $\CC$.
\item The cofibrations of $L_{(H/\VV)}\CC$ are exactly those of $\CC$.
\item The fibrant object of $L_{(H/\VV)}\CC$ are those fibrant $(H/\VV)$-local objects $Z$ of $\CC$.
\item The weak equivalences of $L_{(H/\VV)}\CC$ are the $(H/\VV)$-local equivalences.
\end{enumerate}
\begin{proof} Let $S$ be an $\XX$-small set of cofibrations between cofibrant objects representing all and only the homotopy classes of $H$. Choose a generating set of cofibrations $I$ for $\VV$ with cofibrant domains. Set
\begin{equation*}
L_{(H/\VV)}\CC:=L_{I\Box S}\CC,
\end{equation*}
the left Bousfield localization of $\CC$ by $I\Box S$ (\ref{thm_lbexist}). By \cite[Proposition 17.4.16]{MR2003j:18018}, an object $Z$ that is fibrant in $\CC$ is $(I\Box S)$-local if and only if for any $\fromto{X}{Y}$ in $I$ and any $\fromto{A}{B}$ in $S$, the diagram
\begin{equation*}
\xymatrix@C=18pt@R=18pt{\RMor_{\VV}(Y,\MOR_{\CC}(B,Z))\ar[d]\ar[r]&\RMor_{\VV}(X,\MOR_{\CC}(B,Z))\ar[d]\\
\RMor_{\VV}(Y,\MOR_{\CC}(A,Z))\ar[r]&\RMor_{\VV}(X,\MOR_{\CC}(A,Z))}
\end{equation*}
is a homotopy pullback. Thus $Z$ is $(I\Box S)$-local if and only if for any morphism $\fromto{A}{B}$ of $S$, the induced morphism $\fromto{\MOR_{\CC}(B,Z)}{\MOR_{\CC}(A,Z)}$ is homotopy right orthogonal \cite[Definition 17.8.1]{MR2003j:18018} to every element of $I$. By \cite[Theorem 17.8.18]{MR2003j:18018}, this is equivalent to the condition that $\fromto{\MOR_{\CC}(B,Z)}{\MOR_{\CC}(A,Z)}$ is a weak equivalence in $\VV$. Thus the fibrant objects of $L_{(H/\VV)}\CC$ are exactly those fibrant objects $Z$ of $\CC$ such that the morphism $\fromto{\RMOR_{\CC}(B,Z)}{\RMOR_{\CC}(A,Z)}$ is a weak equivalence in $\VV$ for every $\fromto{A}{B}$ representing an element of $H$.

The inheritance of left properness and the $\XX$-combinatoriality follows from the general theory of left Bousfield localizations (\ref{thm_lbexist}).

The cofibrations are unchanged; hence the the unit axiom holds, and the pushout-product $i\Box f$ of a cofibration $i:\fromto{X}{Y}$ of $\VV$ with a cofibration $f:\fromto{A}{B}$ of $\CC$ is a cofibration of $\CC$ that is trivial if $i$ is. To show that $L_{(H/\VV)}\CC$ is a model $\VV$-category, it thus suffices to show that if $f$ is a trivial cofibration of $L_{(H/\VV)}\CC$, then $i\Box f$ is a weak equivalence. By \ref{lem_Q2varsgencof}, it suffices to verify this for $f$ an element of a generating set of trivial cofibrations of $L_{(H/\VV)}\CC$. By \ref{prp_HoveyJLHMJ}, $L_{(H/\VV)}\CC$ has an $\XX$-small set of generating trivial cofibrations with cofibrant domains; so let $J$ denote such a set, and let $f\in J$. Now by adjunction, one verifies that $i\Box f$ is a weak equivalence if, for any fibrant object $Z$ of $L_{(H/\VV)}\CC$, the diagram
\begin{equation*}
\xymatrix@C=18pt@R=18pt{\RMor_{\VV}(Y,\MOR_{\CC}(B,Z))\ar[d]\ar[r]&\RMor_{\VV}(X,\MOR_{\CC}(B,Z))\ar[d]\\
\RMor_{\VV}(Y,\MOR_{\CC}(A,Z))\ar[r]&\RMor_{\VV}(X,\MOR_{\CC}(A,Z))}
\end{equation*}
is a homotopy pullback. Since both of the vertical morphisms are weak equivalences, the desired result follows.

The right Quillen functor $\fromto{L_{(H/\VV)}\CC}{\CC}$ induces a fully faithful $(\Ho\VV)$-functor
\begin{equation*}
\fromto{\Ho L_{(H/\VV)}\CC}{\Ho \CC}.
\end{equation*}
The characterization of weak equivalences now follows from \ref{lem_RMORdetectsws}.

Now suppose $\DD$ a model $\VV$-category and $F:\fromto{\CC}{\DD}$ a left Quillen $\VV$-functor such that for any $f\in S$, $Ff$ is a weak equivalence. Then $Ff$ is a trivial cofibration of $\DD$, and for any $i\in I$, the morphism $F(i\Box f)=i\Box Ff$ is also a trivial cofibration. Hence any such $F$ factors uniquely through $\fromto{\mathbf{C}}{L_{(H/\mathbf{V})}\mathbf{C}}$ by the universal property enjoyed by ordinary left Bousfield localizations. This completes the proof.
\end{proof}
\end{thm}

\begin{prp}\label{prp_elblocsymmon} Suppose that $\CC$, $\VV$, and $H$ together satisfy each of the conditions (\ref{thm_elbexist}.\ref{CClpropXcomb}) through  (\ref{thm_elbexist}.\ref{VVcofgen}), and, in addition, each of the following conditions.
\begin{enumerate}[(\ref{prp_elblocsymmon}.A)]
\item The tuple $(\CC,\otimes_{\CC}^{\CC},\MOR_{\CC}^{\CC})$ is a symmetric monoidal model category.
\item\label{cofdomIforCC} The model category $\CC$ is tractable.
\item There exists a set (not necessarily $\XX$-small) $G$ of cofibrant homotopy generators of $\CC$, with the property that for any element $A\in G$ and any fibrant $(H/\VV)$-local object $B$ of $\CC$, $\MOR_{\CC}^{\CC}(A,B)$ is $(H/\VV)$-local.
\end{enumerate}
Then the enriched Bousfield localization $L_{(H/\VV)}\CC$ is a symmetric monoidal $\XX$-tractable model $\XX$-category with the symmetric monoidal structure of $\CC$.
\begin{proof} Since the cofibrations are unchanged, $I$ is an $\XX$-small set of generating cofibrations for $L_{(H/\VV)}\CC$ with cofibrant domains. It suffices to verify that for any trivial cofibration $i:\fromto{X}{Y}$ of $L_{(H/\VV)}\CC$ and any element $f:\fromto{A}{B}$ of $I$, the pushout-product $i\Box_{\CC}^{\CC}f$ is a weak equivalence. By \cite[Proposition 17.4.16]{MR2003j:18018}, this holds in turn if, for any fibrant object $Z$ of $L_{(H/\VV)}\CC$, the diagram
\begin{equation*}
\xymatrix@C=18pt@R=18pt{\RMor_{\CC}(Y,\MOR_{\CC}^{\CC}(B,Z))\ar[d]\ar[r]&\RMor_{\CC}(X,\MOR_{\CC}^{\CC}(B,Z))\ar[d]\\
\RMor_{\CC}(Y,\MOR_{\CC}^{\CC}(A,Z))\ar[r]&\RMor_{\CC}(X,\MOR_{\CC}^{\CC}(A,Z))}
\end{equation*}
is a homotopy pullback of $\VV$. The horizontal morphisms are weak equivalences if the objects $\MOR_{\CC}(A,Z)$ and $\MOR_{\CC}(B,Z)$ are $(H/\VV)$-local. This follows from the observation that $A$ and $B$ are homotopy colimits of objects of $G$, and that $(H/\VV)$-locality is preserved under homotopy limits.
\end{proof}
\end{prp}

\subsection*{Application I: Postnikov towers for spectral model categories}{Spectral model categories are model categories enriched in the symmetric monoidal model category of symmetric spectra.}

\begin{ntn} Suppose $\XX$ a universe. Write $\mathbf{Sp}_{\XX}^{\Sigma}$ for the stable symmetric monoidal model category of symmetric spectra in pointed $\XX$-small simplicial sets \cite{}. For any integer $j$ and any symmetric spectrum $E$, write $\pi_jE$ for the $j$-th stable homotopy group of a fibrant replacement of $E$ in $\mathbf{Sp}_{\XX}^{\Sigma}$.
\end{ntn}

\begin{nul} Suppose $\MM$ an $\mathbf{Sp}_{\XX}^{\Sigma}$-model category.
\end{nul}

\begin{dfn} For any integer $n$, an object $X$ of $\MM$ is \emph{$n$-truncated} if for any object $Z$ of $\MM$ and any $j>n$, the group $\pi_j\RMOR_{\MM}(Z,X)=0$.
\end{dfn}

\begin{prp}\label{prp:stntruncated} For any integer $n$, there exists a combinatorial, left proper, $\mathbf{Sp}_{\XX}^{\Sigma}$ model structure on the category $\MM$ --- the \emph{$n$-truncated model structure} $\MM_{\leq n}$ --- satisfying the following conditions.
\begin{enumerate}[(\ref{prp:ntruncated}.1)]
\item The cofibrations of $\MM_{\leq n}$ are precisely the cofibrations of $\MM$.
\item The fibrant objects of $\MM_{\leq n}$ are precisely the fibrant, $n$-truncated objects of $\MM$.
\item The weak equivalences between the fibrant objects are precisely the weak equivalences of $\MM$.
\end{enumerate}
\begin{proof} Let $G$ be an $\XX$-small set of cofibrant homotopy generators of $\MM$. One verifies easily that the $n$-truncated model structure is the enriched left Bousfield localization $L_{(H(n)/\mathbf{Sp}_{\XX}^{\Sigma})}\MM$ with respect to the set
\begin{equation*}
H(n):=\{\fromto{S^j\otimes X}{X}\ |\ X\in G, 0\leq n<j\}.\qedhere
\end{equation*}
\end{proof}
\end{prp}

\begin{ntn} Denote by $\ZZ$ the category whose objects are integers, in which there is a unique morphism $\fromto{m}{n}$ if and only if $m\leq n$.
\end{ntn}

\begin{prp}\label{prp:stpostnikov} There exists a combinatorial, left proper, simplicial model structure on the presheaf category $\MM(\ZZ)$ --- the \emph{stable Postnikov model structure} $\MM(\ZZ)_{\Omega^{-\infty}\Post}$ --- satisfying the following consitions.
\begin{enumerate}[(\ref{prp:stpostnikov}.1)]
\item The cofibrations are the objectwise cofibrations.
\item\label{item:stpostnikovfibrant} The fibrant objects are those sequences
\begin{equation*}
\xymatrix@C=18pt{\dots\ar[r]&X(n)\ar[r]&\dots\ar[r]&X(1)\ar[r]&X(0)\ar[r]&X(-1)\ar[r]&\dots\ar[r]&X(-n)\ar[r]&\dots}
\end{equation*}
satisfying the following conditions.
\begin{enumerate}[(\ref{prp:stpostnikov}.\ref{item:stpostnikovfibrant}.1)]
\item Each object $X(n)$ is fibrant.
\item Each morphism $\fromto{X(n)}{X(n-1)}$ is a fibration.
\item Each object $X(n)$ is $n$-truncated.
\item For any object $Z$ of $\MM$ and any $j<n$, the morphism $\fromto{X\langle n\rangle}{X\langle n-1\rangle}$ induces an isomorphism
\begin{equation*}
\fromto{\pi_j\RMOR_{\MM}(Z,X(n))}{\pi_j\RMOR_{\MM}(Z,X(n-1))}.
\end{equation*}
\end{enumerate}
\item The weak equivalences between the fibrant objects are precisely the objectwise weak equivalences.
\end{enumerate}
\begin{proof} Let $G$ be an $\XX$-small set of cofibrant homotopy generators of $\MM$. As in \ref{prp:postnikov}, the stable Postnikov model structure is the enriched left Bousfield localization $L_{(H/\mathbf{Sp}^{\Sigma})}\MM(\ZZ)_{\inj}$ of the injective model structure with respect to the set
\begin{equation*}
H:=\{\fromto{R^n(S^j\otimes X)}{R^nX}, \fromto{R^{j-1}(S^n\otimes X)}{R^j(S^n\otimes X)}\ |\ X\in G, n<j\}.\qedhere
\end{equation*}
\end{proof}
\end{prp}

\begin{dfn}\label{dfn:stPosttowerhypercomplete}
\begin{enumerate}[(\ref{dfn:stPosttowerhypercomplete}.1)]
\item The \emph{stable Postnikov tower} of an object $X$ in $\MM$ is a fibrant model for
\begin{equation*}
\const X:=[\xymatrix@C=18pt{\dots\ar@{=}[r]&X\ar@{=}[r]&X\ar@{=}[r]&X\ar@{=}[r]&\dots}]
\end{equation*}
in $\MM(\ZZ)_{\Omega^{-\infty}\Post}$; it will be denoted
\begin{equation*}
\xymatrix@C=18pt{\dots\ar[r]&X\langle n\rangle\ar[r]&\dots\ar[r]&X\langle 1\rangle\ar[r]&X\langle 0\rangle\ar[r]&X\langle -1\rangle\ar[r]&\dots\ar[r]&X\langle -n\rangle\ar[r]&\dots}
\end{equation*}
\item One says that a morphism $\fromto{X}{Y}$ of $\MM$ is \emph{$\infty$-connected} if it induces an equivalence $\fromto{\const X}{\const Y}$ in $\MM(\ZZ)_{\Omega^{-\infty}\Post}$.
\item One says that $\MM$ is \emph{stably hypercomplete} if every object is the homotopy limit of its stable Postnikov tower, i.e., if $\holim\circ\LL\const$ is isomorphic to the identity functor on $\Ho\MM$.
\end{enumerate}
\end{dfn}

\begin{prp} There exists a \emph{stable hypercompletion} $\fromto{\MM}{\MM^{\wedge}}$ of $\MM$, which is the initial object in the category of left Quillen $(\mathbf{Sp}_{\XX}^{\Sigma})$-functors $\fromto{\MM}{\NN}$ under which $\infty$-connected morphisms are sent to weak equivalences.
\begin{proof} Simply define $\MM^{\wedge}$ as the enriched left Bousfield localization of $\MM$ with respect to the set $\{\fromto{\lim R\const X}{X}\ |\ X\in G\}$, where $G$ is again an $\XX$-small set of cofibrant homotopy generators of $\MM$, and $R$ is a fibrant replacement functor.
\end{proof}
\end{prp}

\subsection*{Application II: Local model structures}{As a final application, I describe the local model structures on categories of presheaves valued in a symmetric monoidal model category.}

\begin{nul} Suppose $\XX$ a universe, suppose $(C,\tau)$ an $\XX$-small site, and suppose $\VV$ an $\XX$-tractable symmetric monoidal model category with cofibrant unit $\mathbf{1}_{\VV}$.
\end{nul}

\begin{ntn} Write $y:\fromto{C}{\mathbf{Set}_{\XX}(C)}$ for the usual Yoneda embedding, and write $y_{\VV}:\fromto{C}{\VV(C)}$ for the $\VV$-enriched Yoneda embedding, defined by copowers:
\begin{equation*}
y_{\VV}X:=yX\cdot\mathbf{1}_{\VV}:\goesto{Y}{\Mor_C(Y,X)\cdot\mathbf{1}_{\VV}}
\end{equation*}
\end{ntn}

\begin{prp} The category $\VV(C)$, with either its injective or projective model structure, is a $\VV$-model category.
\begin{proof} Suppose $i:\fromto{K}{L}$ a (trivial) cofibration of $\VV$; then for any objectwise (trivial) cofibration (respectively, any objectwise (trivial) fibration) $f:\fromto{X}{Y}$ of $\VV(C)$, the morphism $i\Box f$ (resp., $\mor_{\Box}(i,f)$) is an objectwise (trivial) cofibration (resp., an objectwise (trivial) fibration).
\end{proof}
\end{prp}

\begin{prp} The injective model category $\VV(C)_{\inj}$ is symmetric monoidal.
\begin{proof} Since both cofibrations and weak equivalences are defined objectwise, the pushout-product axiom is immediate.
\end{proof}
\end{prp}

\begin{prp} The projective model category $\VV(C)_{\proj}$ is symmetric monoidal if and only if there exists a set $I$ of generating cofibrations and a set $J$ of generating trivial cofibrations for $\VV$ satisfying the following condition: for any element $[\fromto{X}{Y}]\in I$ (respectively, any element $[\fromto{X}{Y}]\in J$), and any objects $K,L\in\Obj C$, the morphism
\begin{equation*}
\fromto{(\Mor_C(-,K)\times\Mor_C(-,L))\cdot X}{(\Mor_C(-,K)\times\Mor_C(-,L))\cdot Y}
\end{equation*}
is a cofibration (resp. trivial cofibration) of $\VV(C)$.
\begin{proof} Clearly if the projective model structure is symmetric monoidal, any cofibration (resp., trivial cofibration) must satisfy the condition demanded of $I$ (resp., $J$).

Conversely, suppose $I$ and $J$ have been chosen to meet this condition. Then set
\begin{eqnarray}
I_{\VV(\Obj C)}&:=&\bigcup_{K\in\Obj C}(I\times\Prod_{L\neq K}\id_{\varnothing});\nonumber\\
J_{\VV(\Obj C)}&:=&\bigcup_{K\in\Obj C}(J\times\Prod_{L\neq K}\id_{\varnothing}).\nonumber
\end{eqnarray}
Thus $I_{\VV(\Obj C)}$ (resp., $J_{\VV(\Obj C)}$) is a set of generating cofibrations (resp., trivial cofibrations) for $\VV(\Obj C)$, and thus
\begin{equation*}
I_{\VV(C)}:=e_!I_{\VV(\Obj C)}\textrm{\quad (resp.,\quad}J_{\VV(C)}:=e_!J_{\VV(\Obj C)}\textrm{\quad)}
\end{equation*}
is a set of generating cofibrations (resp., trivial cofibrations) for $\VV(C)_{\proj}$. One now easily verifies that the set $I_{\VV(C)}$ (resp., $J_{\VV(C)}$) is the set
\begin{eqnarray}
I_{\VV(C)}&=&\{\fromto{\Mor_C(-,K)\cdot X}{\Mor_C(-,K)\cdot Y}\ |\ K\in\Obj C, [\fromto{X}{Y}]\in I\};\nonumber\\
J_{\VV(C)}&=&\{\fromto{\Mor_C(-,K)\cdot X}{\Mor_C(-,K)\cdot Y}\ |\ K\in\Obj C, [\fromto{X}{Y}]\in J\}.\nonumber
\end{eqnarray}
One thus verifies that the condition of the proposition is precisely the statement that for any element $[i:\fromto{S}{T}]\in I_{\VV(C)}$ (resp., any morphism $[\fromto{S}{T}]\in J_{\VV(C)}$) and any object $L\in\Obj C$, the morphism $\fromto{S\otimes y_{\VV}(L)}{T\otimes y_{\VV}(L)}$ is a cofibration (resp., trivial cofibration), whence it follows from the pushout-product axiom for the $\VV$ enrichment of $\VV(C)_{\proj}$ that for any fibration $p:\fromto{V}{U}$, the morphism
\begin{equation*}
\MOR_{\Box}(i,p):\fromto{\MOR_{\VV(C)}(T,V)}{\MOR_{\VV(C)}(S,V)\times_{\MOR_{\VV(C)}(S,U)}\MOR_{\VV(C)}(T,U)}
\end{equation*}
is a fibration (resp., trivial fibration) that is trivial if $p$ is.
\end{proof}
\end{prp}

\begin{cor} If $C$ has all products, then $\VV(C)_{\proj}$ is a symmetric monoidal model category.
\end{cor}

\begin{dfn} An $\VV$-valued presheaf $F:\fromto{C^{\op}}{\VV}$ is said to \emph{satisfy $\tau$-descent} if for any $\tau$-covering sieve $[\fromto{R}{yX}]\in\tau(X)$, the morphism
\begin{equation*}
\fromto{FX}{\holim_{Y\in(C/R)^{\op}}FY}
\end{equation*}
is an isomorphism of $\Ho\VV$. In this case, $F$ will be called an \emph{$\VV$-valued sheaf.}\footnote{It may be the case that some mathematicians would prefer that I call these ``stacks.'' My insistence on calling them sheaves is meant to be a reflection of the principle that descent for a presheaf is really \emph{always} relative to a homotopy theory on the category in which it takes values.}
\end{dfn}

\begin{nul} It should be noted that the $\VV$-valued sheaves do not satisfy hyperdescent in general; the condition above is only the requirement that a $\VV$-valued sheaf satisfy so-called \emph{\v{C}ech descent}. Ensuring that the $\VV$-valued sheaves satisfy descent with respect to all $\tau$-hypercoverings usually requires a further localization, which I leave for interested readers to formulate in general. In the case of simplicial presheaves, this further localization is simply the hypercompletion, and in the case of spectral presheaves it is the spectral hypercompletion.
\end{nul}

\begin{thm}\label{thm:locprojlocinj} There exist two $\XX$-tractable $\VV$-model structures on the $\VV$-category $\VV(C)$ --- the \emph{$\tau$-local projective} model structure $\VV(C,\tau)_{\proj}$ and (respectively) the \emph{$\tau$-local injective} model structure $\VV(C,\tau)_{\inj}$ --- satisfying the following conditions.
\begin{enumerate}[(\ref{thm:locprojlocinj}.1)]
\item The cofibrations are exactly the projective (resp., injective) cofibrations.
\item The fibrant objects are the projective (resp., injective) fibrant sheaves.
\item The weak equivalences between fibrant objects are precisely the objectwise weak equivalences.
\end{enumerate}
\begin{proof} Set
\begin{equation*}
H:=\{\fromto{R\cdot\mathbf{1}_{\VV}}{y_{\VV}X}\ |\ [\fromto{R}{yX}]\in\tau(X)\}.
\end{equation*}
Then set
\begin{eqnarray}
\VV(C,\tau)_{\proj}&:=&L_{(H/\VV)}\VV(C)_{\proj};\nonumber\\
\VV(C,\tau)_{\inj}&:=&L_{(H/\VV)}\VV(C)_{\inj}.\nonumber
\end{eqnarray}

It now suffices to show that the fibrant objects are as described, and it suffices to do this for the $\tau$-local projective model category $\VV(C,\tau)_{\proj}$. For any $\tau$-covering sieve $[\fromto{R}{yX}]\in\tau(X)$, write $R\cong\colim_{Y\in(C/R)}yY$; one verifies easily that the corresponding colimit $R\cdot\mathbf{1}_{\VV}\cong\colim_{Y\in(C/R)}y_{\VV}Y$ is a homotopy colimit. Thus a fibrant object of $\VV(C,\tau)_{\proj}$ is an objectwise fibrant $\VV$-valued presheaf $F$ such that
\begin{equation*}
FX\simeq\fromto{\RMOR^{\VV}_{\VV(C)_{\proj}}(y_{\VV}X,F)}{\RMOR^{\VV}_{\VV(C)_{\proj}}(R\cdot\mathbf{1}_{\VV},F)}\simeq\holim_{Y\in(C/R)^{\op}}FY
\end{equation*}
for any $\tau$-covering sieve $[\fromto{R}{yX}]\in\tau(X)$.
\end{proof}
\end{thm}

\begin{prp} Suppose $\pi:\fromto{(C,\tau)}{(D,\upsilon)}$ a morphism of sites (hence a functor $\pi^{-1}:\fromto{D}{C}$). Then there is a Quillen adjunction
\begin{equation*}
\pi^{\star}=(\pi^{-1})_!:\xymatrix@1@C=18pt{\VV(D,\upsilon)_{\proj}\ar@<1ex>[r]&\VV(C,\tau)_{\proj}\ar@<1ex>[l]}:(\pi^{-1})^{\star}=\pi_{\star}.
\end{equation*}
\begin{proof} This follows immediately from the universal property of the enriched left Bousfield localization.
\end{proof}
\end{prp}

\begin{thm} The local injective model category $\VV(C,\tau)_{\inj}$ is symmetric monoidal, and if the projective model category $\VV(C)_{\proj}$ is symmetric monoidal, then so is $\VV(C,\tau)_{\proj}$.
\begin{proof} The proofs of the statements are identical, since by \ref{prp_elblocsymmon}, it suffices to show that there exists a set $G$ of cofibrant homotopy generators of $\VV$ such that for any element $Z\in G$, any object $W\in C$, and any local injective fibrant object $Y\in\VV(C)$, the presheaf $\MOR_{\VV(C)}(Z\otimes y_{\VV}W,Y)$ satisfies $\tau$-descent.

To verify this, suppose $X$ an object of $C$, and suppose $[\fromto{R}{yX}]\in\tau(X)$; then one verifies easily that the morphism
\begin{equation*}
\underset{U\in(C/R)^{\op}}{\colim}\fromto{\mathbf{1}_{\VV}\cdot(yW\times yU)}{\mathbf{1}_{\VV}\cdot(yW\times yX)}
\end{equation*}
is a weak equivalence between cofibrant objects. Hence
\begin{equation*}
\fromto{Z\otimes_{\VV(C)}^{\VV}\underset{U\in(C/R)^{\op}}{\colim}\mathbf{1}_{\VV}\cdot(yW\times yU)}{Z\otimes_{\VV(C)}^{\VV}(\mathbf{1}_{\VV}\cdot(yW\times yX))},
\end{equation*}
and therefore
\begin{equation*}
\fromto{\MOR_{\VV(C)}^{\VV}((Z\otimes_{\VV(C)}^{\VV}y_{\VV}W)\otimes y_{\VV}X,Y)}{\underset{U\in(C/R)^{\op}}{\holim}\MOR_{\VV(C)}^{\VV}((Z\otimes_{\VV(C)}^{\VV}y_{\VV}W)\otimes y_{\VV}U,Y)},
\end{equation*}
are weak equivalences of $\VV$, whence the desired descent statement.
\end{proof}
\end{thm}

\bibliographystyle{amsplain}
\bibliography{math}

\end{document}